\documentclass[12pt,oneside,english]{amsart}
\usepackage[latin1]{inputenc}
\usepackage{amssymb}

\makeatletter
\newtheorem{theorem}{Theorem}
\newtheorem{lemma}{Lemma}
\newtheorem{example}{Example}

\newtheorem{remark}{Remark}

\newtheorem{corollary}{Corollary}

\usepackage{babel}

\makeatother
\begin{document}
\baselineskip=17pt

\title[On the basis polynomials]{On the basis polynomials in the theory of
permutations with prescribed up-down structure}

\author{Vladimir Shevelev}
\address{Departments of Mathematics \\Ben-Gurion University of the
 Negev\\Beer-Sheva 84105, Israel. e-mail:shevelev@bgu.ac.il}


\begin{abstract}
Let $\pi=(\pi_1,\pi_2,\hdots,\pi_n)$ be permutation of the elements $1,2,\hdots,n. $
Positive integer $k\leq2^{n-1}$ we call index of $\pi,$ if in its binary notation
as $n$-digital binary number, the 1's correspond to the ascent points.
We study behavior and properties of numbers of permutations of $n$ elements having index $k.$
\end{abstract}

\maketitle

\section{Introduction}

     D.Andre \cite{2} first considered (1881) the problem of the
enumerating the alternating permutations $\pi=(\pi_1,\ldots \pi_n)$
of the numbers $1,2,\ldots, n$ for which ups and downs are
alternating:
$$
\pi_1<\pi_2>\pi_3< \ldots
$$

     This problem has a highly aesthetic solution: the exponential
generating function of such permutations is the sum of tangent and
secant. But only after a century (1968) I.Niven \cite{11} considered
a general problem of the enumerating the permutations with given
up-down structure. For permutation $\pi=(\pi_1,\ldots,\pi_n)$, the
sequence $(q_1,q_2,\ldots,q_{n-1})$, where
\begin{equation}\label{1}
q_j=sign(\pi_{j+1}-\pi_j)=\begin{cases}1,& if \;\; \pi_{j+1}>\pi_j
\\-1,& if\;\; \pi_{j+1}<\pi_j\end{cases},
\end{equation}

is called a \slshape Niven's signature.\upshape For example,
$a=(2,1,5,4,3)$ has the signature $(-1,1,-1,-1)$.

     Denote $[q_1,q_2,\ldots,q_{n-1}]$ the number of permutations
having the Niven's signature $(q_1,q_2,\ldots,q_{n-1})$. In view of
the symmetry, we have

\begin{equation}\label{2}
[q_1,q_2,\ldots,q_{n-1}]=[-q_{n-1},-q_{n-2},\ldots,-q_1].
\end{equation}

Niven obtained the following basic result.
\newpage

\begin{theorem}\label{t1} \cite{11}. Let in the signature
$(q_1,q_2,\ldots,q_{n-1})$ the indices of those $q_i$ which are $+1$
be $k_1<k_2<\ldots<k_m$ (if such $q_i$ do not exist then assume
$m=0$). Put in addition $k_0=0,\;\;k_{m+1}=n$. Then

\begin{equation}\label{3}
[q_1,q_2,\ldots,q_{n-1}]= det N,
\end{equation}

where $N=\{n_{ij}\}$   is the square matrix of order $m+1$ in which

\begin{equation}\label{4}
n_{ij}=\begin{pmatrix} k_i\\k_{j-1}\end{pmatrix},\;\;i,j=1,2,\ldots,
m+1.
\end{equation}
\end{theorem}

     After this celebrated Niven's result and until now there has
been a series of articles by many authors. We mention only eleven
papers in chronological order: N.G.Bruijn, 1970 \cite{5},
H.O.Foulkes, 1976 \cite{7}, L.Carlitz, 1978 \cite{6}, G.Viennot,
1979 \cite{18}, C.L.Mallows and L.A.Shepp, 1985 \cite{9}, V.Arnold,
1990 \cite{3}, V.S.Shevelev, 1996 \cite{14}, G.Szpiro, 2001
\cite{17},B.Shapiro, M.Shapiro and A.Vainshtein, 2005 \cite{12},
F.C.S.Brown, T.M.A. Fink and K.Willbrand, 2007 \cite{4},R. Stanley, 2007 \cite{16}.

     According to the de Bruijn-Viennot algorithm \cite{5},
\cite{18}, the calculation of $[q_1,q_2,\ldots,q_{n-1}]$ could be
realized using the following numerical triangle. At the top of the
triangle put 1. Write 0 to the right (left) if $q_{n-1}$ is $1(-1)$.
For example, if $q_{n-1}=1$ then the first two elements of the
triangle are: $\begin{matrix}& 1
\\ && 0 \end{matrix}$. Summing these elements we write the sum to the
left:      $\begin{matrix}& 1
\\1 && 0 \end{matrix}$. Now the following $0$ we write to the right
of the last element (to the left of the first element) if $q_{n-2}$
is $1(-1)$. For example, if $q_{n-2}=1$ then we have
$\begin{matrix}& 1
\\1 && 0 \\&&& 0\end{matrix}$. Now the  third row is obtained by summing
each element in the third row from the right to the left, first
element being $0$, with the elements left and above it in the second
row:
$$
\begin{matrix}&& 1
\\&1 && 0 \\1 && 0 && 0 \end{matrix}
$$

\; \;\;\;\;

In case of $q_{n-2}=-1$ we have $\begin{matrix}&& 1
\\&1 && 0 \\0  \end{matrix}$ and the third row is obtained by
summing each element in the third row from the left to the right, the
first element being $0$, with the elements right and above it in the
second row:
\newpage
$$
\begin{matrix}&& 1
\\&1 && 0 \\0 && 1 && 1 \end{matrix}
$$

     The process is continued until the n-th row which corresponds to
$q_1$. Now the sum of elements of the n-th row is equal to
$[q_1,q_2,\ldots,q_{n-1}]$.

\begin{example}.  For signature $(-1,1,1,-1,1)$ we have the triangle

$$
\setcounter{MaxMatrixCols}{21}\begin{matrix}&&&&& 1
\\&&&&1 && 0 \\&&&0 && 1 && 1
\\&&2&&2&&1&&0\\&5&&3&&1&&0&&0\\0&&5&&8&&9&&9&&9\end{matrix}
$$
\;\;\;
\end{example}

     Therefore, $[-1,1,1,-1,1]=0+5+8+9+9+9=40$.
In order to obtain a \slshape weight generalization, \upshape let us
consider a matrix function which we call "alternant" (cf.\cite{13}).
If a permutation $\pi$ has the signature $(q_1,q_2,\ldots,q_{n-1}),$
then we write $\pi\in (q_1,q_2,\ldots,q_{n-1})$. Furthermore, if
$\pi_i=j,$ then to the two-dimensional point $(i,j)$ assign the
"weight" $a_{ij}.$

Let $A=(a_{ij})$ be an $n\times n$ matrix. Denote

\begin{equation}\label{5}
alt_{(q_1,\ldots,q_{n-1})}A=\sum_{\pi\in(q_1,\ldots,q_{n-1})}\prod^n_{i=1}a_{i\pi_i}.
\end{equation}

Let $A_{1j},\; j=1,\ldots,n,$ be $(n-1)\times(n-1)$ matrix which is
obtained from $A$ by the deletion of the first row and the j-th
column. Denote $A^{(+1)}_{1j}(A^{(-1)}_{1j})$ the matrix which is
obtained from $A_{1j}$ by replacing the $j-1$ first (the n-j last)
elements of its first row by 0's. Then, by the Viennot's algorithm,
we deduce the following expansion of the alternant by the first row
of the matrix.

\begin{theorem}\label{t2}  $($cf.\enskip\cite{13}$)$

\begin{equation}\label{6}
alt_{(q_1,\ldots,q_{n-1})}A=\sum^n_{j=1}a_{1j}alt_{(q_2,\ldots,q_{n-1})}A^{(q_1)}_{1j}.
\end{equation}
\end{theorem}

Note that, if $A=J_n-n\times n$ matrix composed of 1's only, then

\begin{equation}\label{7}
alt_{(q_1,\ldots,q_{n-1})}A=[q_1,\ldots,q_{n-1}].
\end{equation}\newpage

In case of arbitrary (0,1) matrix $A$, Theorem 2 gives enumeration
the permutations having signature $(q_1,\ldots,q_{n-1})$ with
restriction on positions. For example, if $I_n$ is $(n\times n)$
identity matrix, then $alt_{(q_1,\ldots,q_{n-1})}(J-I)$ gives the
number of such permutations without fixed points.

Note that, by (\ref{5})

\begin{equation}\label{8}
alt_{(1)}\begin{pmatrix}a_{11}&a_{12}\\a_{21}&a_{22}\end{pmatrix}=a_{11}a_{22},\;\;
alt_{(-1)}\begin{pmatrix}a_{11}&a_{12}\\a_{21}&a_{22}\end{pmatrix}=a_{12}a_{21}.
\end{equation}

\begin{example}\label{2}.
$$
alt_{(1,-1,1)}J_4=alt_{(1,-1,1)}\begin{pmatrix}1&1&1&1\\1&1&1&1\\1&1&1&1\\
1&1&1&1\end{pmatrix}=alt_{(-1,1)}\begin{pmatrix}1&1&1\\1&1&1\\1&1&1
\end{pmatrix}+
$$

$$
alt_{(-1,1)}\begin{pmatrix}0&1&1\\1&1&1\\1&1&1
\end{pmatrix}+alt_{(-1,1)}\begin{pmatrix}0&0&1\\1&1&1\\1&1&1
\end{pmatrix}+alt_{(-1,1)}\begin{pmatrix}0&0&0\\1&1&1\\1&1&1
\end{pmatrix}=
$$

$$
\left(alt_{(1)}\begin{pmatrix}0&0\\1&1
\end{pmatrix}+alt_{(1)}\begin{pmatrix}1&0\\1&1
\end{pmatrix}+alt_{(1)}\begin{pmatrix}1&1\\1&1
\end{pmatrix}\right)+
$$

$$
\left(alt_{(1)}\begin{pmatrix}1&0\\1&1
\end{pmatrix}+alt_{(1)}\begin{pmatrix}1&1\\1&1
\end{pmatrix}\right)+alt_{(1)}\begin{pmatrix}1&1\\1&1
\end{pmatrix}=0+1+1+1+1+1=5,
$$

while

$$
alt_{(1,-1,1)}(J_4-I)=alt_{(1,-1,1)}\begin{pmatrix}0&1&1&1\\1&0&1&1\\1&1&0&1\\
1&1&1&0\end{pmatrix}=alt_{(-1,1)}\begin{pmatrix}0&1&1\\1&0&1\\1&1&0
\end{pmatrix}+
$$

$$
alt_{(-1,1)}\begin{pmatrix}0&0&1\\1&1&1\\1&1&0
\end{pmatrix}+alt_{(-1,1)}\begin{pmatrix}0&0&0\\1&1&0\\1&1&1
\end{pmatrix}=
$$
$$
\left(alt_{(1)}\begin{pmatrix}1&0\\1&0
\end{pmatrix}+alt_{(1)}\begin{pmatrix}1&0\\1&1
\end{pmatrix}\right)+alt_{(1)}\begin{pmatrix}1&1\\1&1
\end{pmatrix}=2.
$$
\end{example}
Note that alternant is also useful for enumeration the permutations
with some additional conditions. For example, if it is necessary to
enumerate the permutations $\pi$ with signature
$(q_1,q_2,\ldots,q_{n-1})$ for which $\pi_1=l,\;\;\pi_n=m$ then we
should calculate $alt J_n^{(l,m)}$ where $J_n^{(l,m)}$ is obtained
from $J_n$ by replacing all 1's of the first and the last rows by
0's except the l-th 1 and \newpage the m-th 1 correspondingly. For example,
there are only 2 such permutations in the case of the signature
$(-1,1,1,-1,1),\;n=6,\;l=2,\;m=6$.

Let us now introduce an \slshape index \upshape of the Niven's
signature in the following way: the integer $k=k_n$ is called the
index of the signature $(q_1,q_2,\ldots,q_{n-1})$ if $(n-1)$-digit
binary representation of $k$ is

\begin{equation}\label{9}
k=\sum^{n-1}_{i=1}q_i^{\shortmid} 2^{n-i-1},
\end{equation}

where

\begin{equation}\label{10}
q^\shortmid_i=\begin{cases}
1,\;\;if\;\;q_i=1,\\0,\;\;otherwise\end{cases}.
\end{equation}

(cf \cite{13}, \cite{14}, where we accepted $k+1$ as an index). Denote
$S^{(k)}_n$ the set of permutations of elements $1,2,\ldots,n$
having the index $k$, and put

\begin{equation}\label{11}
\left\{\begin{matrix} n \\ k
\end{matrix}\right\}=\left|S_n^{(k)}\right\}.
\end{equation}

Let $k\in[2^{t-1}, 2^t)$ and the $(n-1)$-digit binary expansion of
$k$ (\ref{9}) has a form:

\begin{equation}\label{12}
k=\underbrace{0 \ldots 0}_{n-t-1}\; 1 \;\underbrace{0 \ldots
0}_{s_2-s_1-1}\; 1 \;\underbrace{0 \ldots 0}_{s_3-s_2-1}\; 1\ldots 1
\;\underbrace{0 \ldots 0}_{s_m-s_{m-1}-1}\; 1 \;\underbrace{0 \ldots
0}_{t-s_m}\enskip,
\end{equation}

where

\begin{equation}\label{13}
1=s_1< s_2 <\ldots < s_m
\end{equation}
are places of 1's \slshape after \upshape $n-t-1$ 0's before the
first 1.

     From Theorem 1 follows an important formula.

\begin{theorem}\label{t3}
For index $k\in [2^{t-1}, 2^t)\enskip(12),$ we have
\end{theorem}

\begin{equation}\label{14}
\left\{\begin{matrix} n \\ k
\end{matrix}\right\}=
\end{equation}\newpage
$$\begin{vmatrix}1 & 1& 0&  \hdots &
0\\1&\begin{pmatrix}n-t+s_2-1 \\ s_2 -1 \end{pmatrix} & 1&
\hdots & 0 \\1&\begin{pmatrix}n-t+s_3-1 \\ s_3 -1 \end{pmatrix} &
\begin{pmatrix}n-t+s_3-1 \\ s_3 -s_2 \end{pmatrix}&  \hdots &
0\\\hdots & \hdots & \hdots & \hdots & \hdots  \\
1 & \begin{pmatrix}n-t+s_m-1 \\ s_m -1 \end{pmatrix}&
\begin{pmatrix}n-t+s_m-1 \\ s_m -s_2 \end{pmatrix} & \begin{pmatrix}n-t+s_m-1 \\
s_m -s_3 \end{pmatrix} & 1\\
1 & \begin{pmatrix}n \\ t \end{pmatrix} & \begin{pmatrix}n \\
t+1-s_2
\end{pmatrix}& \begin{pmatrix}n \\ t+1-s_3 \end{pmatrix} & \begin{pmatrix}n \\
t+1-s_m \end{pmatrix}\end{vmatrix}.$$

In \cite{14}, using a techniques of permanents, the following explicit
formula was proved.

\begin{theorem}\label{t4}  \cite {14}.  For $k\in[2^{t-1}, 2^t)$ we
have

$$\left\{\begin{matrix} n \\ k
\end{matrix}\right\}=(-1)^m+$$
\begin{equation}\label{15}
\sum^m_{p=1}(-1)^{m-p}\sum_{1\leq i,< i_2< \hdots< i_p\leq m}
\left( \begin{matrix} n \\ t+1-s_{i_p}\end{matrix}\right)\prod
^p_{r=2}\begin{pmatrix} n-t+s_{i_r}-1
\\ s_{i_r}-s_{i_{r-1}}\end{pmatrix}.
\end{equation}
\end{theorem}

Most likely, (\ref{15}) is the first non-determinant formula in a
closed form for the number of permutations with prescribed up-down
structure (cf.\cite{16}, \cite{4}).

The Thue-Morse sequence \cite{10}, \cite{8} is defined by

\begin{equation}\label{16}
\tau_n=(-1)^{\sigma(n)},
\end{equation}

where $\sigma(n)$ denotes the number of 1's in the binary
expansion of $n$.

     Thus, from (\ref{15}) immediately follows an interesting
arithmetical property of $\left\{\begin{matrix} n \\ k
\end{matrix}\right\}$.

\begin{theorem}\label{t5} \cite{14}. If all different from 1
divisors of $n$ are larger than $\lfloor\log_2 2k\rfloor, $ then
\end{theorem}

\begin{equation}\label{17}
\left\{\begin{matrix} n \\ k
\end{matrix}\right\}\equiv \tau_k \;(\mod{n}).
\end{equation}

\begin{remark}\label{1}
It is evident that the validity of $(\ref{17})$ does not depend on the
fact whether $k$ is a constant or a function of $n$.
\end{remark}\newpage

\begin{remark}\label{2}
From Theorem $5$ and $(\ref{16})$ it follows that, if $n$ has only
sufficiently large different from 1 divisors, then the number of
permutations of $n$ elements with arbitrary prescribed up-down
structure is $\pm 1(\mod{n})$. Recently, a special case of this
result for a prime $n$ was reproduced in \cite{4}.
\end{remark}

\section{Basis polynomials}

As it follows from (\ref{15}), if $k$ does not depend on $n,$ then
$\left\{\begin{matrix}n\\ k \end{matrix}\right\}$ is a polynomial in
$n$ of the degree

\begin{equation}\label{18}
t=\lfloor\log_2(2k)\rfloor.
\end{equation}

Indeed, the degree of the polynomial in the interior sum of
(\ref{15})equals to

$$
t+1-s_{i_p}+(s_{i_2}-s_{i_1})+(s_{i_3}-s_{i_2})+\hdots+(s_{i_p}-s_{i_{p-1}})
=t+1-s_{i_1}\leq t
$$

and the equality is attained in those summands of the sum in which
$i_1=1$.

     Let us draw an analogy with the binomial coefficients
$\begin{pmatrix}n\\ k \end{pmatrix}.$

1a. $\begin{pmatrix}n\\ k \end{pmatrix}$ is the number of subsets of
the cardinality $k$ of a set of $n$ elements.

1b. $\left\{\begin{matrix} n\\k \end{matrix}\right \}$  is the
number of permutations of $n$ elements having the up-down index $k$.

2a.  Each subset of a set of $n$ elements is contained in the number
of $\begin{pmatrix}n\\ k \end{pmatrix}$ subsets for some value of
$k$.

2b.  Each permutation of $n$ elements is contained in the number of
 $\left\{\begin{matrix}n\\ k
\end{matrix}\right\}$  permutations for some value of the up-down index $k$.

3a.  $\sum^n_{k=0}\left(\begin{matrix}n\\ k
\end{matrix}\right)=2^n$.

3b. $ \sum^{2^{n-1}-1}_{k=0}\left\{\begin{matrix}n\\ k
\end{matrix}\right\}=n! $

4a. $\begin{pmatrix}n\\ n-k \end{pmatrix}=\begin{pmatrix}n\\ k
\end{pmatrix}$.

4b.  In view of (\ref{2}),  $\left\{\begin{matrix}n\\ k
\end{matrix}\right\}=\left\{\begin{matrix}n\\ 2^{n-1}-1-k
\end{matrix}\right\}$.

5a.$\begin{pmatrix}n\\ 0 \end{pmatrix}=\begin{pmatrix}n\\ n
\end{pmatrix}=1$.\newpage

5b.  $\left\{\begin{matrix}n\\ 0
\end{matrix}\right\}=\left\{\begin{matrix}n\\ 2^{n-1}-1
\end{matrix}\right\}=1$.

The latter equality corresponds to the identity permutation.

6.  The central binomial coefficients and the "central" numbers
$\left\{\begin{matrix}n\\ k
\end{matrix}\right\}$ are equal one to another. Indeed, below (section 3) we prove that

$$
\left\{\begin{matrix}2n\\2^{n-1}-1
\end{matrix}\right\}=\begin{pmatrix}2n-1\\n-1
\end{pmatrix},
$$

$$
\left\{\begin{matrix}2n+1\\2^n-1
\end{matrix}\right\}=\begin{pmatrix}2n\\n
\end{pmatrix}.
$$

In view of this analogy, we call $\left\{\begin{matrix}n\\k
\end{matrix}\right\},\;\;k=1,2,\hdots$, the \slshape basis polynomial \upshape in
theory of permutations with prescribed up-down structure.

     Note that $\left\{\begin{matrix}n\\k
\end{matrix}\right\}$, just as $\begin{pmatrix}n\\k \end{pmatrix}$,
generally is not a polynomial if $k$ is a function of $n$.

\begin{example}. In case of alternating permutations
$\pi_1<\pi_2>\pi_3<\hdots, $ we have the sequence of indices
$\{k_{n-1}\}$ such that

$$
k_1=1,\;\;k_2=2,\;\;k_3=5,\;\;k_4=10,\;\;k_5=21,\hdots\enskip.
$$
\end{example}

Here

$$
k_n-k_{n-2}=2^{n-1},\;\; n\geq 3,
$$
whence

\begin{equation}\label{19}
k_{n-1}=\frac{2^{n+1}-3+(-1)^{n}}{6},\;\;n=1,2,\hdots\enskip
\end{equation}

Thus, from the classical Andre's result we obtain

\begin{equation}\label{20}
\sum^{\infty}_{n=0}\left\{\begin{matrix}n\\\frac{2^{n+1}-3+(-1)^{n}}{6}
\end{matrix}\right\}\frac{x^n}{n!}=\tan{x}+\sec{x},
\end{equation}

where we put $\left\{\begin{matrix}0\\0
\end{matrix}\right\}=1.$

From (\ref{20}) we have some values of $\left\{\begin{matrix}n\\m
\end{matrix}\right\}:$

$$\left\{\begin{matrix}1\\0
\end{matrix}\right\}=1,\;\left\{\begin{matrix}2\\1
\end{matrix}\right\}=1,\;\left\{\begin{matrix}3\\2
\end{matrix}\right\}=2,\;$$\newpage
\begin{equation}\label{21}
\left\{\begin{matrix}4\\5
\end{matrix}\right\}=5,\;\left\{\begin{matrix}5\\10
\end{matrix}\right\}=16,\;\left\{\begin{matrix}6\\21
\end{matrix}\right\}=61,\;\hdots\enskip .
\end{equation}

It is well-known that these values are explicitly expressed with
help of the absolute values of the Bernoulli and Euler numbers.

\section{Another general formula for basis polynomials}

Here we use formulas (\ref{12})-(\ref{14}) for obtaining more simple
general explicit formula. It is important,for the subsequent
development of our theory, to present the results in the form of
linear combination of the binomial coefficients $\begin{pmatrix}n\\i
\end{pmatrix}$. Firstly consider several special cases.

a) $m=1,\; k=2^{t-1}$. Then

\begin{equation}\label{22}
\left\{\begin{matrix}n\\2^{t-1}
\end{matrix}\right\}=\left|\begin{matrix}1& 1\\1 &\begin{pmatrix}n\\t
\end{pmatrix}
\end{matrix}\right|=\begin{pmatrix}n\\t
\end{pmatrix}-1;
\end{equation}

b) $m=2,\; k=2^{t-1}+2^{u-1}$. Here $s_1=1,\;s_2=t-u+1$ and we have

$$
\left\{\begin{matrix}n\\2^{t-1}+2^{u-1}
\end{matrix}\right\}=\left|\begin{matrix}1& 1 & 0\\1
&\begin{pmatrix}n-u\\t-u
\end{pmatrix}& 1\\1 &\begin{pmatrix}n\\t
\end{pmatrix}& \begin{pmatrix}n\\u
\end{pmatrix}\end{matrix}\right|=\begin{pmatrix}n-u\\t-u
\end{pmatrix}\begin{pmatrix}n\\u
\end{pmatrix}-\begin{pmatrix}n\\t
\end{pmatrix}-\begin{pmatrix}n\\u
\end{pmatrix}+1.
$$

Notice that,

\begin{equation}\label{23}
\begin{pmatrix}n-u\\t-u
\end{pmatrix}\begin{pmatrix}n\\u
\end{pmatrix}=\begin{pmatrix}t\\u
\end{pmatrix}\begin{pmatrix}n\\t
\end{pmatrix}.
\end{equation}

Therefore,

\begin{equation}\label{24}
\left\{\begin{matrix}n\\2^{t-1}+2^{u-1}
\end{matrix}\right\}=\left(\begin{pmatrix}t\\u
\end{pmatrix}-1\right)\begin{pmatrix}n\\t
\end{pmatrix}-\begin{pmatrix}n\\u
\end{pmatrix}+1.
\end{equation}

c)  $m=3,\;k=2^{t-1}+2^{u-1}+2^{v-1}$. Here $s_1=1,\; s_2=t-u+1,\;
s_3=t-v+1$ and we have\newpage

$$
\left\{\begin{matrix}n\\2^{t-1}+2^{u-1}+2^{v-1}
\end{matrix}\right\}=\left|\begin{matrix}1 & 1 & 0 & 0\\1 &
\begin{pmatrix} n-u\\ t-u \end{pmatrix}& 1 & 0\\1 &
\begin{pmatrix} n-v\\ t-v \end{pmatrix}& \begin{pmatrix} n-v\\ u-v \end{pmatrix} &
1\\1 & \begin{pmatrix} n\\t \end{pmatrix} & \begin{pmatrix} n\\u
\end{pmatrix} &
\begin{pmatrix} n\\v \end{pmatrix}
\end{matrix}\right|=
$$

$$
=\left|\begin{matrix} \begin{pmatrix} n-u\\ t-u \end{pmatrix} & 1 &
0\\\begin{pmatrix} n-v\\ t-v \end{pmatrix}& \begin{pmatrix} n-v\\
u-v\end{pmatrix}& 1\\\begin{pmatrix} n\\ t \end{pmatrix}&
\begin{pmatrix} n\\u \end{pmatrix}&\begin{pmatrix} n\\ v
\end{pmatrix}\end{matrix}\right|-\left|\begin{matrix} 1 & 1 &
0\\1& \begin{pmatrix} n-v\\
u-v\end{pmatrix}& 1\\1&
\begin{pmatrix} n\\u \end{pmatrix}&\begin{pmatrix} n\\ v
\end{pmatrix}\end{matrix}\right|=
$$

$$
=\begin{pmatrix} n-u\\
t-u\end{pmatrix}\begin{pmatrix} n-v\\
u-v\end{pmatrix}\begin{pmatrix} n\\
v\end{pmatrix}-\begin{pmatrix} n-u\\
t-u\end{pmatrix}\begin{pmatrix} n\\
u\end{pmatrix}-\begin{pmatrix} n-v\\
t-v\end{pmatrix}\begin{pmatrix} n\\
v\end{pmatrix}-
$$

\begin{equation}\label{25}
-\begin{pmatrix} n-v\\
u-v\end{pmatrix}\begin{pmatrix} n\\
v\end{pmatrix}+\begin{pmatrix} n\\
t\end{pmatrix}+\begin{pmatrix} n\\
u\end{pmatrix}+\begin{pmatrix} n\\
v\end{pmatrix}-1.
\end{equation}

Notice that,

\begin{equation}\label{26}
\begin{pmatrix} n-u\\
t-u\end{pmatrix}\begin{pmatrix} n-v\\
u-v\end{pmatrix}\begin{pmatrix} n\\
v\end{pmatrix}=\begin{pmatrix} t\\
u\end{pmatrix}\begin{pmatrix} u\\
v\end{pmatrix}\begin{pmatrix} n\\
t\end{pmatrix},
\end{equation}

$$
\begin{pmatrix} n-u\\
t-u\end{pmatrix}\begin{pmatrix} n\\
u\end{pmatrix}=\begin{pmatrix} t\\
u\end{pmatrix}\begin{pmatrix} n\\
t\end{pmatrix},
$$

\begin{equation}\label{27}
\begin{pmatrix} n-v\\
t-v\end{pmatrix}\begin{pmatrix} n\\
v\end{pmatrix}=\begin{pmatrix} t\\
v\end{pmatrix}\begin{pmatrix} n\\
t\end{pmatrix},
\end{equation}

$$
\begin{pmatrix} n-v\\
u-v\end{pmatrix}\begin{pmatrix} n\\
v\end{pmatrix}=\begin{pmatrix} u\\
v\end{pmatrix}\begin{pmatrix} n\\
u\end{pmatrix}.
$$

Therefore,

$$
\left\{\begin{matrix} n\\
2^{t-1}+2^{u-1}+2^{v-1}\end{matrix}\right\}=\left(\begin{pmatrix} t\\
u\end{pmatrix}\begin{pmatrix} u\\
v\end{pmatrix}-\begin{pmatrix} t\\
u\end{pmatrix}-\begin{pmatrix} t\\
v\end{pmatrix}+1\right)\begin{pmatrix} n\\
t\end{pmatrix}-
$$

\begin{equation}\label{28}
-\left(\begin{pmatrix} u\\
v\end{pmatrix}-1\right)\begin{pmatrix} n\\
u\end{pmatrix}+\begin{pmatrix} n\\
v\end{pmatrix}-1.
\end{equation}

At last, for arbitrary $m$, we obtain the following theorem.\newpage

\begin{theorem}\label{t6}
If
\begin{equation}\label{29}
k=2^{t_1-1}+2^{t_2-1}+\hdots+2^{t_m-1},\;\;t_1>t_2>\hdots>t_m,
\end{equation}

then

$$
\left\{\begin{matrix}
n\\k\end{matrix}\right\}=(-1)^m\left(1-\sum^m_{i=1}\begin{pmatrix}
n\\t_i \end{pmatrix}+\sum_{1\leq i<j\leq
m}\begin{pmatrix}n\\t_i\end{pmatrix}\begin{pmatrix}t_i\\t_j\end{pmatrix}-\right.
$$

\begin{equation}\label{30}
-\sum_{1\leq i<j\leq l\leq
m}\begin{pmatrix}n\\t_i\end{pmatrix}\begin{pmatrix}t_i\\t_j\end{pmatrix}\begin{pmatrix}t_j\\t_l\end{pmatrix}
+\hdots\left.
+(-1)^m\begin{pmatrix}n\\t_1\end{pmatrix}\prod^{m-1}_{j=1}\begin{pmatrix}t_i\\t_{i+1}
\end{pmatrix}\right).
\end{equation}
\end{theorem}
\bfseries Proof.\mdseries \enskip Theorem 6 directly follows from (12), (15), (29) and the following easily proved identity (if to take into account that from (12) and (29) we have $s_i=t+1-t_i,\enskip i=1,....m):$
$$\begin{pmatrix} n\\t_{i_p}\end{pmatrix} \prod^p_{r=2}\begin{pmatrix} n-t_{i_r}\\ t_{i_{r-1}}-t_{i_r}\end{pmatrix}=
\begin{pmatrix} n\\ t_{i_1}\end{pmatrix}\begin{pmatrix} t_{i_1}\\ t_{i_2}\end{pmatrix}
\begin{pmatrix} t_{i_2}\\ t_{i_3}\end{pmatrix}\hdots\begin{pmatrix} t_{i_{p-1}}\\ t_{i_p}\end{pmatrix}.$$ $\blacksquare$

It is clear that in (\ref{30}) $\left\{\begin{matrix}
n\\k\end{matrix}\right\}$ is presented as a linear combination of
$\begin{pmatrix} n\\i \end{pmatrix},\;\;0\leq i\leq t_1$. More
exactly, as (\ref{30}) shows, $\left\{\begin{matrix}
n\\k\end{matrix}\right\}$ is a linear combination of elements of the
last row of the determinant (\ref{14}), that does not follow from
(\ref{14}) directly.

Indeed, in case of (\ref{29}) we have

$$
s_1=1,\;\;\; s_i=t_1-t_i+1,\;\;\; i=2,3,\hdots, m.
$$

In (\ref{14}) $t_1=t$ and, consequently,

$$
t+1-s_i=t_i,\;\;\;i=1,2,\hdots, m.
$$

Thus, $\left\{\begin{matrix} n\\k \end{matrix}\right\}$ is an
alternating sum of some elementary symmetric polynomials of binomial
coefficients. It is also a polynomial in $n$ of degree
$t_1=\lfloor\log_2(2k)\rfloor$.

Note that another form of Theorem 6 is

\begin{theorem} \label{t7}   In conditions $(\ref{29})$ we have

\begin{equation}\label{31}
\left\{\begin{matrix} n\\k
\end{matrix}\right\}=(-1)^m+\sum^m_{p=1}c_p\begin{pmatrix} n\\t_p
\end{pmatrix},
\end{equation}
\newpage
where
$$
c_p=(-1)^m\left(-1+\sum^m_{j=p+1}\begin{pmatrix} t_p\\t_j
\end{pmatrix} -\sum_{p+1\leq j<l\leq m}\begin{pmatrix} t_p\\t_j
\end{pmatrix}\begin{pmatrix} t_j\\t_l
\end{pmatrix}+ \hdots +\right .
$$
\begin{equation}\label{32}
\left. +(-1)^{m-p-1}\prod^m_{j=p+1}\begin{pmatrix} t_{j-1}\\t_j
\end{pmatrix}\right).
\end{equation}
\end{theorem}
In particular,
$$
c_m=(-1)^{m+1},
$$
\begin{equation}\label{33}
c_{m-1}=(-1)^m\left(-1+\begin{pmatrix} t_{m-1}\\ t_m
\end{pmatrix}\right),
\end{equation}
$$
c_{m-2}=(-1)^m\left(-1+\begin{pmatrix} t_{m-2}\\ t_{m-1}
\end{pmatrix}+\begin{pmatrix} t_{m-2}\\ t_m
\end{pmatrix}-\begin{pmatrix} t_{m-2}\\ t_{m-1}
\end{pmatrix}\begin{pmatrix} t_{m-1}\\ t_m
\end{pmatrix}\right),
$$

etc.

\begin{example}\label{4}

Let $k=2^m-1$. Then, according to $(\ref{29}),$ we have

$$
t_m=1,\;\;\; t_{m-1}=2,\;\;\hdots,\;\;\; t_1=m
$$

and, by $(\ref{33}),$ we find
$$
c_m=(-1)^{m+1},
$$
$$
c_{m-1}=(-1)^m(-1+2)=(-1)^m,
$$
$$
c_{m-2}=(-1)^m(-1+3+3-3\cdot 2)=(-1)^{m+1},
$$

and, by induction,
$$
c_m=-c_{m-1}= c_{m-2}=\hdots=(-1)^{m-1}c_1=(-1)^{m-1},
$$
i.e.
$$
c_p=(-1)^{p-1}.
$$
\end{example}

Thus, by (\ref{31}) we have

\begin{equation}\label{34}
\left\{\begin{matrix} n\\ 2^m
-1\end{matrix}\right\}=(-1)^m+\sum^m_{p=1}(-1)^{p-1}\begin{pmatrix}
n\\m-p+1 \end{pmatrix}=\sum^m_{j=0}(-1)^{m-j}\begin{pmatrix} n\\j
\end{pmatrix}=\begin{pmatrix} n-1\\m
\end{pmatrix}.
\end{equation}

The latter identity is proved easily by induction over $m$.

     In particular, putting in (\ref{34}) $n=2m$ and $n=2m+1$ we
have

\begin{equation}\label{35}
\left\{\begin{matrix} 2m\\ 2^m-1\end{matrix}\right\}=\begin{pmatrix}
2m-1\\m\end{pmatrix},\;\;\left\{\begin{matrix} 2m+1\\
2^m-1\end{matrix}\right\}=\begin{pmatrix} 2m\\m\end{pmatrix}.
\end{equation}\newpage

This proves the analogy for the "central" number
$\left\{\begin{matrix} n\\ k\end{matrix}\right\}$ and the central
binomial coefficients.

Comparing with (\ref{14}), we obtain an identity:

\begin{equation}\label{36}
\begin{vmatrix}1 & 1& 0& 0&\hdots & 0
\\1&\begin{pmatrix}n-m+1 \\ 1 \end{pmatrix} & 1& 0 &
\hdots & 0 \\1&\begin{pmatrix}n-m+2 \\ 2 \end{pmatrix} &
\begin{pmatrix}n-m+2 \\ 1 \end{pmatrix}& 1 & \hdots &
0\\\hdots & \hdots & \hdots & \hdots & \hdots & \hdots \\
1 & \begin{pmatrix}n-1 \\ m-1 \end{pmatrix}&
\begin{pmatrix}n-1 \\ m-2 \end{pmatrix} & \begin{pmatrix}n-1 \\
m-3\end{pmatrix} & \hdots & 1\\
1 &  \begin{pmatrix}n \\m\end{pmatrix}& \begin{pmatrix}n \\ m-1
\end{pmatrix}
& \begin{pmatrix}n \\
m-2 \end{pmatrix}& \hdots & \begin{pmatrix}n \\ 1
\end{pmatrix}\end{vmatrix}=\begin{pmatrix}n-1 \\ m \end{pmatrix}
\end{equation}

     Formulas (\ref{31})-(\ref{33}) alow to calculate rather effectively basis
polynomials $\left\{\begin{matrix} n\\ k\end{matrix}\right\}$.
Nevertheless, there exist a recursions which are more
effective for calculations.

\section{Recursion relation for basis polynomials}

Consider now $\left\{\begin{matrix} n\\ k\end{matrix}\right\}$ from
a more formal point of view as a polynomial (\ref{30}) or
(\ref{31}). If $k$ has larger digits than $n-1$ then
$\left\{\begin{matrix} n\\ k\end{matrix}\right\}$ loses its
combinatorial sense and could take even negative values.
Nevertheless, the formal values of $\left\{\begin{matrix} n\\
k\end{matrix}\right\}$ are useful since the coefficients $c_p$
(\ref{32}) could be represented as some values of the basis
polynomials.

\begin{theorem}\label{t8}
In conditions $(\ref{29})$ in $(\ref{31})$ we have

$$
c_p=\left\{\begin{matrix} t_p\\ k-2^{t_p-1}\end{matrix}\right\}.
$$
\end{theorem}

\bfseries Proof.\mdseries  \enskip By (\ref{29}),

$$
k_1=k-2^{t_p-1}=2^{t_1-1}+2^{t_2-1}+\hdots +
$$
$$
+2^{t_{p-1}-1}+2^{t_{p+1}-1}+\hdots + 2^{t_m-1}.
$$

Using (\ref{30}) for $k_1$ and substituting $n=t_p\enskip,$ we obtain\newpage

$$
\left\{\begin{matrix} t_p\\
k-2^{t_p-1}\end{matrix}\right\}=(-1)^{m-1}\left(1-\sum^m_{i=p+1}\begin{pmatrix}
t_p\\ t_i \end{pmatrix}+\sum_{p+1\leq i < j \leq m}\begin{pmatrix} t_p\\
t_i\end{pmatrix}\begin{pmatrix} t_i\\
t_j\end{pmatrix}- \hdots \right)
$$

and the comparison with (\ref{32}) gives the theorem. $\blacksquare$

From Theorems 7, 8 we obtain a very simple recursion relation.

\begin{theorem}\label{t9}    If $k=2^{t_1-1}+2^{t_2-1}+ \hdots +
2^{t_m-1},$ then

\begin{equation}\label{37}
\left\{\begin{matrix} n\\
k\end{matrix}\right\}=(-1)^m+\sum^m_{p=1}\left\{\begin{matrix} t_p\\
k-2^{t_p-1}\end{matrix}\right\}\begin{pmatrix}n\\
t_p\end{pmatrix}.
\end{equation}
\end{theorem}

\begin{example}     Knowing $\left\{\begin{matrix} n\\
j\end{matrix}\right\},\;\; j\leq 20,$  to find $\left\{\begin{matrix} n\\
21\end{matrix}\right\}$.  We have

$$
21= 2^{5-1}+2^{3-1}+ 2^{1-1},\;\; t_1=5,\;\; t_2=3, \;\;t_3=1.
$$
\end{example}

By (\ref{37}), we obtain

\begin{equation}\label{38}
\left\{\begin{matrix} n\\
21\end{matrix}\right\}=-1+\left\{\begin{matrix} 5\\
5\end{matrix}\right\}\begin{pmatrix}n\\
5\end{pmatrix}+\left\{\begin{matrix} 3\\
17\end{matrix}\right\}\begin{pmatrix}n\\
3\end{pmatrix}+\left\{\begin{matrix} 1\\
20\end{matrix}\right\}\begin{pmatrix}n\\
1\end{pmatrix}.
\end{equation}

Using formulas (see Appendix)

$$
\left\{\begin{matrix} n\\
5\end{matrix}\right\}=2\begin{pmatrix}n\\
3\end{pmatrix}-\begin{pmatrix}n\\
1\end{pmatrix}+1,
$$
$$
\left\{\begin{matrix} n\\
17\end{matrix}\right\}=4\begin{pmatrix}n\\
5\end{pmatrix}-\begin{pmatrix}n\\
1\end{pmatrix}+1,
$$
$$
\left\{\begin{matrix} n\\
20\end{matrix}\right\}=9\begin{pmatrix}n\\
5\end{pmatrix}-\begin{pmatrix}n\\
3\end{pmatrix}+1,
$$

we conclude that
$$
\left\{\begin{matrix} n\\
5\end{matrix}\right\}=2\cdot 10 -5 +1=16,
$$
$$
\left\{\begin{matrix} 3\\
17\end{matrix}\right\}=-2,
$$
$$
\left\{\begin{matrix} 1\\
20\end{matrix}\right\}=1
$$

and, by (\ref{38}), we find

$$
\left\{\begin{matrix} n\\
21\end{matrix}\right\}=16\begin{pmatrix}n\\
5\end{pmatrix}-2\begin{pmatrix}n\\
3\end{pmatrix}+\begin{pmatrix}n\\
1\end{pmatrix}-1.
$$
\newpage
\section{Another determinant formula for basis polynomials}

In conditions (\ref{29}) the determinant (\ref{14}) has the form
(cf.3,a),b),c)):

\begin{equation}\label{39}
\left\{\begin{matrix} n\\
k\end{matrix}\right\}=\begin{vmatrix} 1& 1& 0 & 0& \hdots & 0\\
1& \begin{pmatrix}n-t_2\\t_1-t_2\end{pmatrix} & 1& 0& \hdots & 0\\
1& \begin{pmatrix}n-t_3\\t_1-t_3\end{pmatrix}& \begin{pmatrix}n-t_3\\t_2-t_3\end{pmatrix}&
 1& \hdots & 0\\ \hdots & \hdots & \hdots & \hdots & \hdots & \hdots
 \\
1& \begin{pmatrix}n-t_m\\t_1-t_m\end{pmatrix}& \begin{pmatrix}n-t_m\\t_2-t_m\end{pmatrix}&
\begin{pmatrix}n-t_m\\t_3-t_m\end{pmatrix}&  \hdots & 1\\
1 &\begin{pmatrix}n\\
t_1\end{pmatrix}&\begin{pmatrix}n\\
t_2\end{pmatrix}&\begin{pmatrix}n\\
t_3\end{pmatrix}& \hdots &\begin{pmatrix}n\\
t_m\end{pmatrix}\end{vmatrix}.
\end{equation}

Note that this determinant possesses an astonishing property. If to
replace the lower triangular submatrix with the main diagonal of 1's
by the upper one such that the elements $\begin{pmatrix} n-t_j\\
t_i-t_j\end{pmatrix}\;(i< j)$ are mapped to elements $\begin{pmatrix} t_i\\
t_j\end{pmatrix}$ which are symmetric respectively the diagonal of
1's (and which do not depend on $n(!)$), then the determinant does
not change its value. If, in addition, to interchange the places of
the first and last rows, then we obtain the following result.

\begin{theorem}\label{t10}
$k=2^{t_1-1}+2^{t_2-1}+ \hdots + 2^{t_m-1},\;\;t_1>t_2>\hdots>t_m$,
then
\end{theorem}

\begin{equation}\label{40}
\left\{\begin{matrix} n\\
k\end{matrix}\right\}=(-1)^m \begin{vmatrix} 1& \begin{pmatrix}n\\t_1\end{pmatrix}&
\begin{pmatrix}n\\t_2\end{pmatrix}& \begin{pmatrix}n\\t_3\end{pmatrix}& \hdots &
\begin{pmatrix}n\\t_m\end{pmatrix}\\
1& 1 &\begin{pmatrix}t_1\\t_2\end{pmatrix}&
\begin{pmatrix}t_1\\t_3\end{pmatrix}&
\hdots &\begin{pmatrix}t_1\\t_m\end{pmatrix}\\
1& 0& 1 &\begin{pmatrix}t_2\\t_3\end{pmatrix}&\hdots
&\begin{pmatrix}t_2\\t_m\end{pmatrix}\\1 & 0 & 0 & 1 & \hdots &
\begin{pmatrix}t_3\\t_m\end{pmatrix}\\\hdots &\hdots &\hdots &\hdots &\hdots &\hdots
\\1 & 0 & 0 & 0 & \hdots & 1
\end{vmatrix}.
\end{equation}

\bfseries Proof. \mdseries  One can prove this formula using the
analysis of the structure of the diagonals and the comparison with
(\ref{30}). $\blacksquare$

In the case of alternating permutations when $k=k_{n-1}$ is defined by
(\ref{19}), we obtain an $(m+1)\times (m+1)$ determinant
representation of numbers (\ref{21}).\newpage Thus, for $n=2m,$ we have an
identity for the Euler numbers $E_{2m},\;m\geq 1$ (cf.\cite{1},
Table 23.2, \cite{15}, A000364). We drop $(-1)^m$ in order to take
account of the sign of $E_{2m}$.

\begin{equation}\label{41}
E_{2m}=\begin{vmatrix} 1& \begin{pmatrix}2m\\2m-1\end{pmatrix} &
\begin{pmatrix}2m\\2m-3\end{pmatrix}&
\begin{pmatrix}2m\\2m-5\end{pmatrix}& \hdots &
\begin{pmatrix}2m\\1\end{pmatrix}\\
1 & 1 & \begin{pmatrix}2m-1\\2m-3\end{pmatrix}&
\begin{pmatrix}2m-1\\2m-5\end{pmatrix}& \hdots &
\begin{pmatrix}2m-1\\1\end{pmatrix}\\
1 & 0 & 1 & \begin{pmatrix}2m-3\\2m-5\end{pmatrix}& \hdots
&\begin{pmatrix}2m-3\\1\end{pmatrix}\\
1 & 0 & 0 & 1 & \hdots & \begin{pmatrix}2m-5\\1\end{pmatrix}\\
\hdots &\hdots &\hdots &\hdots &\hdots &\hdots \\
1 & 0 & 0 & 0 & \hdots & 1 \end{vmatrix}
\end{equation}

Analogously, putting $n=2m-1,\;\;m\geq 2,$ for the Bernoulli numbers
$B_{2m}$ (cf.\cite{1}, Table 23.2) we have the following determinant
of an $m\times m$ matrix:

\begin{equation}\label{42}
\frac{B_{2m}}{2m}(2^{2m}-1)2^{2m}=\begin{vmatrix} 1&
\begin{pmatrix}2m-1\\2m-2\end{pmatrix} &
\begin{pmatrix}2m-1\\2m-4\end{pmatrix}&\begin{pmatrix}2m-1\\2m-6\end{pmatrix} &
\hdots &\begin{pmatrix}2m-1\\2\end{pmatrix}\\
1 & 1 & \begin{pmatrix}2m-2\\2m-4\end{pmatrix}&
\begin{pmatrix}2m-2\\2m-6\end{pmatrix}& \hdots &
\begin{pmatrix}2m-2\\2\end{pmatrix}\\
1 & 0 & 1 & \begin{pmatrix}2m-4\\2m-6\end{pmatrix}& \hdots
&\begin{pmatrix}2m-4\\2\end{pmatrix}\\
1&0&0&1&\hdots &\begin{pmatrix}2m-6\\2\end{pmatrix}\\
\hdots &\hdots &\hdots &\hdots &\hdots &\hdots \\
1 & 0 & 0 & 0 & \hdots & 1 \end{vmatrix}
\end{equation}

The numbers on the left hand side are the tangent numbers
(\cite{15}, A 000182).

\begin{example} For $m=1,$ we have
$$
E_2=\begin{vmatrix} 1 & \begin{pmatrix}2\\1\end{pmatrix}\\
1 & 1\end{vmatrix}=-1,
$$
for $m=2,$ we have

$$
\frac{B_4}{4}16 \cdot
15=\begin{vmatrix}1 & \begin{pmatrix}3\\2\end{pmatrix}\\
1 & 1\end{vmatrix}=-2,
$$\newpage
which corresponds to $B_4=-\frac{1}{30}$.

     For $m=2,$ we have also
$$
E_4=\begin{vmatrix} 1 & \begin{pmatrix}4\\3\end{pmatrix}&
\begin{pmatrix}4\\1\end{pmatrix}\\
1 & 1& \begin{pmatrix}3\\1\end{pmatrix}\\
1 & 0 & 1 \end{vmatrix}=5,
$$
for $m=3,$ we have
$$
\frac{B_6}{6}64\cdot 63=\begin{vmatrix} 1 &
\begin{pmatrix}5\\4\end{pmatrix}&
\begin{pmatrix}5\\2\end{pmatrix}\\
1 & 1& \begin{pmatrix}4\\2\end{pmatrix}\\
1 & 0 & 1 \end{vmatrix}=16,
$$
which corresponds to $B_6=\frac{1}{42}$.
\end{example}

\section{An identity for partial sums of the basis polynomials}

\begin{theorem}\label{t11}  For $1\leq r\leq n-1,$ we have
\begin{equation}\label{43}
\sum^{2^r-1}_{k=0}\left\{\begin{matrix} n\\ k
\end{matrix}\right\}=n(n-1) \hdots (n-r+1).
\end{equation}
\end{theorem}

\bfseries Proof. \mdseries Sum (\ref{43}) enumerates the
permutations with the $n-r-1$ fixed down points: $1,2,\hdots,
n-r-1$. Let us form an arbitrary permutation $\pi$ of such kind. We
start with position $n-r+1$. We can choose value of $\pi_{n-r-1}$ by
$n$ ways, $\pi_{n-r+2}$ by $n-1$ ways, $\hdots, \pi_n$ by
$n-(r-1)$ways. After that $\pi_1>\pi_2>\hdots > \pi_{n-r}$ are
defined uniquely. Thus, we obtain (\ref{43}).$\blacksquare$

\begin{example}  For $r=3,$ we have (cf Appendix)
\end{example}
$$
\sum^7_{k=0}\left\{\begin{matrix} n\\k
\end{matrix}\right\}=1+\left(\begin{pmatrix} n\\ 1
\end{pmatrix}-1\right)+\left(\begin{pmatrix} n\\ 2
\end{pmatrix}-1\right)+\left(\begin{pmatrix} n\\ 2
\end{pmatrix}-\begin{pmatrix} n\\ 1
\end{pmatrix}+1\right)+
$$
$$
+\left(\begin{pmatrix} n\\ 3
\end{pmatrix}-1\right)+\left(2\begin{pmatrix} n\\ 3
\end{pmatrix}-\begin{pmatrix} n\\ 1
\end{pmatrix}+1\right)+\left(2\begin{pmatrix} n\\ 3
\end{pmatrix}-\begin{pmatrix} n\\ 2
\end{pmatrix}+1\right)+
$$
$$
+\left(\begin{pmatrix} n\\3 \end{pmatrix}-\begin{pmatrix} n\\
2\end{pmatrix}+\begin{pmatrix} n\\ 1
\end{pmatrix}-1\right)=6\begin{pmatrix} n\\3 \end{pmatrix}=n(n-1)(n-2).
$$
\newpage
\section{On the positive integer zeros of the basis polynomials}

\begin{theorem}\label{t12}
If $k=2^{t_1-1}+2^{t_2-1}+\hdots +2^{t_m-1}, \; t_1>t_2>\hdots
>t_m\geq 1,$ then the integers $t_1, t_2, \hdots, t_m$ are roots of
the basis polynomial $\left\{\begin{matrix} n\\
k\end{matrix}\right\}.$
\end{theorem}

\bfseries Proof. \mdseries Substituting in (\ref{40}),
$n=t_j,\;j=1,2,\hdots,m,$ we obtain a determinant with two the same
rows. $\blacksquare$

     Furthermore, the following results are obtained directly from
(\ref{30}):

$$
\left\{\begin{matrix} 1\\
k\end{matrix}\right\}=(-1)^m\left(1-\begin{pmatrix} 1\\t_m
\end{pmatrix}\right)=\begin{cases} (-1)^m,\;\; if \;\; t_m\geq
2,\\0,\;\;\;\; if \;\;\;\; t_m=1\end{cases};
$$

$$
\left\{\begin{matrix} 2\\
k\end{matrix}\right\}=(-1)^m\left(1-\begin{pmatrix} 2\\t_{m-1}
\end{pmatrix}-\begin{pmatrix} 2\\t_m
\end{pmatrix}+\begin{pmatrix} 2\\t_{m-1}
\end{pmatrix}\begin{pmatrix} 2\\t_m
\end{pmatrix}\right)=
$$
$$
=\begin{cases} (-1)^m,\;\;if\;\;t_m\geq 3\\
0,\;\; if\;\;t_m=2\\
0,\;\; if\;\;t_m=1,\;\;t_{m-1}=2,\\
(-1)^{m-1},\;\; if \;\; t_m=1,\;\; t_{m-1}\geq 3\end{cases};
$$

$$
\left\{\begin{matrix} 3\\
k\end{matrix}\right\}=(-1)^m\left(1-\begin{pmatrix}
3\\t_{m-2}\end{pmatrix}-\begin{pmatrix}
3\\t_{m-1}\end{pmatrix}-\begin{pmatrix}
3\\t_m\end{pmatrix}+\begin{pmatrix}
3\\t_{m-2}\end{pmatrix}\begin{pmatrix}
3\\t_{m-1}\end{pmatrix}\right.+
$$
$$
\left.+\begin{pmatrix} 3\\t_{m-2}\end{pmatrix}\begin{pmatrix}
3\\t_m\end{pmatrix}+\begin{pmatrix}
3\\t_{m-1}\end{pmatrix}\begin{pmatrix}
3\\t_m\end{pmatrix}-\begin{pmatrix}
3\\t_{m-2}\end{pmatrix}\begin{pmatrix}
3\\t_{m-1}\end{pmatrix}\begin{pmatrix} 3\\t_m\end{pmatrix}\right)=
$$
$$
=\begin{cases} (-1)^m,\;\; if \;\; t_m\geq 4\\
0, \;\; if \;\; t_m=3\\
0, \;\; if \;\; t_m=2,\;\;t_{m-1}=3\\
2(-1)^{m-1}, \;\; if \;\; t_m=2,\;\;t_{m-1}> 3\\
0, \;\; if \;\; t_m=1,\;\;t_{m-1}=2,\;\;t_{m-2}=3\\
4(-1)^m, \;\; if \;\; t_m=1,\;\;t_{m-1}=2,\;\ t_{m-2}> 3\\
0, \;\; if \;\; t_m=1,\;\;t_{m-1}=3,\;\;t_{m-2}>3\\
2(-1)^{m-1}, \;\; if \;\; t_m=1,\;\;t_{m-1}> 3\end{cases} ;
$$

etc.

Researching the appearance of zero values, we obtain the following result.

\begin{theorem}\label{t13}  Let $ a\in\mathbb{N}.$   If

\begin{equation}\label{44}
k\equiv 2^{a-1}+j\mod{2^a},
\end{equation}

where $j\in[0, 2^{a-1}), $ then $\left\{\begin{matrix} a\\
k\end{matrix}\right\}=0$.
\end{theorem}

\bfseries Proof. \mdseries By (\ref{44}), we have\newpage

$$
k=l 2^a + 2^{a-1} + j,\;\;l \geq 0,\;\; 0 \leq j \leq 2^{a-1}-1.
$$

Therefore, if $k=2^{t_1-1}+2^{t_2-1}+\hdots+2^{t_m},$ then there
exists $p\in [1,m]$ such that $a=t_p$ and the result follows from
Theorem 12. $\blacksquare$

\begin{remark}\label{t14}   Note that, the conversion of Theorem $13$ implies the conversion of Theorem $12.$ This will be obtained at the end of this article (Theorem $22$).
\end{remark}

As a corollary from Theorem 13, it follows a more attractive
statement.

\begin{theorem}\label{t14}     Let
$k=2^{t_1-1}+2^{t_2-1}+\hdots+2^{t_m-1}$ and $1\leq
i\leq\log_2(2k).$ If $i\neq t_p,\;\;p=1,2,\hdots,m,$ then
$\left\{\begin{matrix} i\\ k-2^{i-1}\end{matrix}\right\}=0.$
\end{theorem}

\bfseries Proof. \mdseries   Let $t_l< i < t_{l-1}.$ Then

$$
k-2^{i-1}=2^{t_1-1}+2^{t_2-1}+\hdots+\left(2^{t_{l-1}-1}-2^{i-1}\right)
+2^{t_p-1}+\hdots+2^{t_m-1}=
$$
$$
=2^{t_1-1}+2^{t_2-1}+\hdots+\left(2^{t_{l-1}-2}+2^{t_{l-1}-3}+\hdots+
2^i+2^{i-1}\right)+2^{t_l-1}+\hdots+2^{t_m-1}\equiv
$$
$$
\equiv 2^{i-1}+2^{t_l-1}+\hdots+2^{t_m-1}\pmod{2^i}
$$

and the theorem directly follows from Theorem 13.$\blacksquare$

\section{Another algorithm of evaluation of basis polynomials}

\begin{theorem}\label{t15}   If  $k=2^{t_1-1}+2^{t_2-1}+\hdots+2^{t_m-1},\;
 t_1> t_2>\hdots > t_m\geq 1,$ then

\begin{equation}\label{45}
\left\{\begin{matrix} n\\ k \end{matrix}\right\}=a_1\begin{pmatrix}
n\\ t_m \end{pmatrix}+a_2\begin{pmatrix} n\\ t_{m-1}
\end{pmatrix}+\hdots + a_m\begin{pmatrix}
n\\ t_1 \end{pmatrix}+(-1)^m,
\end{equation}
where integers $a_i,\;\; i=1,2,\hdots, m$, are defined by the system
of the linear equations
\end{theorem}

\begin{equation}\label{46}
\begin{cases} a_1+(-1)^m=0\\
\begin{pmatrix} t_{m-1}\\ t_m \end{pmatrix} a_1+a_2+(-1)^m=0\\
\begin{pmatrix} t_{m-2}\\ t_m \end{pmatrix} a_1+\begin{pmatrix} t_{m-2}\\
 t_{m-1} \end{pmatrix} a_2+a_3+(-1)^m=0\\
 \hdots\hdots\hdots\hdots\hdots\hdots\\
\begin{pmatrix} t_1 \\ t_m \end{pmatrix}a_1+\begin{pmatrix} t_1\\ t_{m-1}
\end{pmatrix}a_2+\hdots +\begin{pmatrix} t_1\\ t_2
\end{pmatrix}a_{m-1}+a_m+(-1)^m=0 \end{cases}.
\end{equation}\newpage

\bfseries Proof. \mdseries  From (\ref{40}) follows a representation
(\ref{45}). Substituting in (\ref{45}) $n=t_m,t_{m-1},\hdots,t_1$
and using Theorem 12, we obtain system (\ref{46}).$\blacksquare$

\begin{example}    Let us find $\left\{\begin{matrix} n\\26
\end{matrix}\right\}.$ We have

$$
26=2^{5-1}+2^{4-1}+2^{2-1}.
$$
\end{example}
Thus, $t_1=5,\;\;t_2=4,\;\;t_3=2,\;\;m=3$. By (\ref{46})

$$
\begin{cases} a_1-1=0\\6a_1+a_2-1=0\\10a_1+5a_2+a_3-1=0\end{cases},
$$
whence $a_1=1,\;\; a_2=-5, \;\; a_3=16.$ Consequently, by (\ref{45})

$$
\left\{\begin{matrix} n\\ 26 \end{matrix}\right\}=\begin{pmatrix}
n\\ 2\end{pmatrix}-5\begin{pmatrix} n\\
4\end{pmatrix}+16\begin{pmatrix} n\\ 5\end{pmatrix}-1.
$$

\section{Another recursion relation for basis polynomials}

\begin{theorem}\label{t16}      Let $k=2^{t_1-1}+2^{t_2-1}+\hdots+2^{t_m-1},\;
 t_1> t_2>\hdots > t_m\geq 1$ . Then, for $l>t_1$ we have
 \end{theorem}

\begin{equation}\label{47}
\left\{\begin{matrix} n\\k+2^{l-1}
\end{matrix}\right\}=\left\{\begin{matrix} l\\k
\end{matrix}\right\}\begin{pmatrix}
n\\ l\end{pmatrix}-\left\{\begin{matrix} n\\k
\end{matrix}\right\}.
\end{equation}

\bfseries Proof. \mdseries  By the latter theorem,

\begin{equation}\label{48}
\left\{\begin{matrix} n\\ k+2^{l-1}\end{matrix}\right\}=
b_1\begin{pmatrix} n\\ t_m\end{pmatrix}+b_2\begin{pmatrix} n\\
t_{m-1}
\end{pmatrix}+\hdots+b_m\begin{pmatrix} n\\ t_1\end{pmatrix}+b_{m+1}\begin{pmatrix} n\\ l\end{pmatrix}
+(-1)^{m+1},
\end{equation}

where $b_i,\;\;i=1,2,\hdots, m+1$, are defined by the following
system:
\begin{equation}\label{49}
\begin{cases} b_1+(-1)^{m+1}=0\\
\begin{pmatrix}t_{m-1}\\ t_m\end{pmatrix} b_1+b_2+(-1)^{m+1}=0\\
\hdots\hdots\hdots\hdots\hdots\hdots\hdots\\
\begin{pmatrix}t_1\\ t_m\end{pmatrix}b_1+\begin{pmatrix}t_1\\
t_{m-1}\end{pmatrix}b_2+\hdots +\begin{pmatrix}t_1\\
t_2\end{pmatrix}b_{m-1} + b_m+ (-1)^{m+1}=0\\
\begin{pmatrix}l\\ t_m\end{pmatrix}b_1+\begin{pmatrix}l\\
t_{m-1}\end{pmatrix}b_2+\hdots +\begin{pmatrix}l\\
t_1\end{pmatrix}b_m + b_{m+1}+ (-1)^{m+1}=0
\end{cases}.
\end{equation}

By comparison of the first $m$ equations of (\ref{49}) with (\ref{46}),
we conclude that\newpage

\begin{equation}\label{50}
b_i=-a_i,\;\; i=1,2,\hdots, m
\end{equation}

and, by the $(m+1)-th$ equation of (\ref{49}) and by  (\ref{45}), we
find that

\begin{equation}\label{51}
b_{m+1}=-\left\{\begin{matrix} l\\k
\end{matrix}\right\}.
\end{equation}

Now from (\ref{45}) , (\ref{48}), (\ref{50}) and (\ref{51}) we
obtain

$$
\left\{\begin{matrix} n\\k
\end{matrix}\right\}+\left\{\begin{matrix} n\\k+2^{l-1}
\end{matrix}\right\}=\left\{\begin{matrix} l\\k
\end{matrix}\right\}\begin{pmatrix} n\\l
\end{pmatrix}
$$
and (\ref{47}) follows. $\blacksquare$

\begin{example}  Starting with  $\left\{\begin{matrix} n\\0
\end{matrix}\right\}=1$ and putting $k=0,\;\;l=1,$ we obtain

$$
\left\{\begin{matrix} n\\1\end{matrix}\right\}=
\left\{\begin{matrix} 1\\0\end{matrix}\right\}\begin{pmatrix} n\\1
\end{pmatrix}-\left\{\begin{matrix} n\\0
\end{matrix}\right\}=\begin{pmatrix} n\\1
\end{pmatrix}-1.
$$
\end{example}
Furthermore, we consecutively find:

$$
putting\;\;k=0,\;\; l=2,\;\left\{\begin{matrix} n\\2
\end{matrix}\right\}=\left\{\begin{matrix} 2\\0
\end{matrix}\right\}\begin{pmatrix} n\\2
\end{pmatrix}-\left\{\begin{matrix} n\\0
\end{matrix}\right\}=\begin{pmatrix} n\\2
\end{pmatrix}-1,
$$
$$
putting\;\;k=1,\;\; l=2,\;\left\{\begin{matrix} n\\3
\end{matrix}\right\}=\left\{\begin{matrix} 2\\1
\end{matrix}\right\}\begin{pmatrix} n\\2
\end{pmatrix}-\left\{\begin{matrix} n\\1
\end{matrix}\right\}=\begin{pmatrix} n\\2
\end{pmatrix}-\begin{pmatrix} n\\1
\end{pmatrix}+1,
$$
$$
putting\;\;k=0,\;\; l=3,\;\left\{\begin{matrix} n\\4
\end{matrix}\right\}=\left\{\begin{matrix} 3\\0
\end{matrix}\right\}\begin{pmatrix} n\\3
\end{pmatrix}-\left\{\begin{matrix} n\\0
\end{matrix}\right\}=\begin{pmatrix} n\\3
\end{pmatrix}-1,
$$
$$
putting\;\;k=1,\;\; l=3,\;\left\{\begin{matrix} n\\5
\end{matrix}\right\}=\left\{\begin{matrix} 3\\1
\end{matrix}\right\}\begin{pmatrix} n\\3
\end{pmatrix}-\left\{\begin{matrix} n\\1
\end{matrix}\right\}=2\begin{pmatrix} n\\3
\end{pmatrix}-\begin{pmatrix} n\\1
\end{pmatrix}+1
$$
etc.

\section{Characteristic conditions for a basis polynomial}

Let $P(n)$ be a polynomial. It is evident that the condition

$$
P(n)=C\begin{pmatrix} n\\k
\end{pmatrix}
$$
with a constant $C$ satisfies if and only if $P(r)=0,\;\;r=0,1,\hdots,k$, and
$k=deg P(n)$. Concerning $\left\{\begin{matrix} n\\k
\end{matrix}\right\}$ we have the following result.  Put

$$
\Delta P(n)=P(n)-P(n-1)
$$

and let $\Delta^r P(n)$ be the $r-th$ difference of $P(n)$.\newpage

\begin{theorem}\label{t17}
For a polynomial $P( n)$ there exists a nonnegative integer $k$ and
a constant $C\neq 0$ such that

\begin{equation}\label{52}
P(n)=C\left\{\begin{matrix} n\\k
\end{matrix}\right\}
\end{equation}

if and only if the following conditions satisfy:

\begin{equation}\label{53}
P(0)\neq 0,\;\;(\Delta^rP(r))P(r)=0,\;\;r=1,2,\hdots,l,
\end{equation}
where $l=deg P(n)$.
\end{theorem}
\bfseries Proof. \mdseries   Note that, we have

\begin{equation}\label{54}
P(n)=P(0)+\sum^l_{r=1}\Delta^rP(r)\begin{pmatrix} n\\r \end{pmatrix}.
\end{equation}
Indeed, put
$$ P(n)=a_0+\sum_{r=1}^la_r\begin{pmatrix} n\\r \end{pmatrix}.$$
Then we consecutively find
$$a_0=P(0),$$
$$ \Delta P(n)=\sum_{r=1}^la_r\begin{pmatrix} {n-1}\\{r-1} \end{pmatrix} , \enskip a_1=\Delta P(1),...,$$
i.e.
$$ \Delta^t P(n)=\sum_{r=1}^la_r\begin{pmatrix} {n-t}\\{r-t} \end{pmatrix} , \enskip a_t=\Delta^t P(t),\enskip t=0,...,l.$$

If all $\Delta^rP(r)=0,\;\;r=1,\hdots,l$, then we put $k=0,\;\;C=P(0)$. If
$\Delta^rP(r)\neq 0,$ for $r=t_1>t_2>\hdots>t_m\geq 1,$ then by (\ref{54})

$$
P(n)=P(0)+\sum^m_{i=1}b_{m+1-i}\begin{pmatrix} n\\t_i \end{pmatrix},
$$
where
$$
b_i=\Delta^{t_i}P(t_i),\;\; i=1,2,\hdots,m.
$$
Putting
$$
\frac{b_i}{P(0)}(-1)^m=a_i,\;\;i=1,2,\hdots,m,
$$
we have

\begin{equation}\label{55}
\frac{(-1)^m}{P(0)}P(n)=(-1)^m+\sum^m_{i=1}a_{m+1-i}\begin{pmatrix}
n\\t_i \end{pmatrix}
\end{equation}\newpage

and, according to (\ref{53}), $P(t_i)=0,\;\;i=1,2,\hdots, m$. Thus, by
Theorem 15, the polynomial (\ref{55}) is $\left\{\begin{matrix} n\\k
\end{matrix}\right\}$ with

$$
k=2^{t_1-1}+2^{t_2-1}+\hdots+2^{t_m-1}.
$$

The converse statement is evident as well, according to Theorem
12.$\blacksquare$

\begin{example}

Consider $P(n)=n^3-3n^2+2n-6$. We have
\end{example}
$$
\Delta P(n)=3n^2-9n+6, \;\;\Delta^2P(n)=6n-12,\;\; \Delta^3P(n)=6,
$$
and see that

$$
(\Delta P(1))P(1)=(\Delta^2P(2))P(2)=(\Delta^3P(3))P(3)=0.
$$
By Theorem 17, we conclude that $P(n)=C\left\{\begin{matrix} n\\k
\end{matrix}\right\}$. Now we easily find $C.$ Since the smallest $r$ for which $ \Delta^r P(3)\neq 0, $ then
$t_1=3$ and $m=1.$  Therefore, $k=4$ and $C=-P(0)=6$. Thus $P(n)=6\left
\{\begin{matrix} n\\4 \end{matrix}\right \}.$

\section{On generating function of the basis polynomials}
Note that, formula (47) gives a possibility to add to $k$ any powers of 2 more than $2^{t_1}.$ Therefore, using some iterations of (47), one can formally to get any $k_1>k.$ Thus we have a natural way to define $\left\{\begin{matrix} n\\k\end{matrix}\right\}$ for $k>2^{n-1}-1.$ Moreover, the following lemma shows that the series $\sum^\infty_{k=0}\left \{\begin{matrix} n\\k\end{matrix}\right \}x^k$ converges in interval $|x|<1$ for every $n.$
\begin{lemma}
For a fixed $n,$ the sequence $\{\left\{\begin{matrix} n\\
k\end{matrix}\right\}\}_{k\geq0}$ is bounded.
\end{lemma}
\bfseries Proof. \mdseries Directly we consecutively have from (\ref{30}):

$$
\left\{\begin{matrix} 1\\
k\end{matrix}\right\}=(-1)^m\left(1-\begin{pmatrix} 1\\t_m
\end{pmatrix}\right)=\begin{cases} (-1)^m,\;\; if \;\; t_m\geq
2,\\0,\;\;\;\; if \;\;\;\; t_m=1\end{cases};
$$

$$
\left\{\begin{matrix} 2\\
k\end{matrix}\right\}=(-1)^m\left(1-\begin{pmatrix} 2\\t_{m-1}
\end{pmatrix}-\begin{pmatrix} 2\\t_m
\end{pmatrix}+\begin{pmatrix} 2\\t_{m-1}
\end{pmatrix}\begin{pmatrix} 2\\t_m
\end{pmatrix}\right)=
$$
$$
=\begin{cases} (-1)^m,\;\;if\;\;t_m\geq 3\\
0,\;\; if\;\;t_m=2\\
0,\;\; if\;\;t_m=1,\;\;t_{m-1}=2,\\
(-1)^{m-1},\;\; if \;\; t_m=1,\;\; t_{m-1}\geq 3\end{cases};
$$
\newpage

$$
\left\{\begin{matrix} 3\\
k\end{matrix}\right\}=(-1)^m\left(1-\begin{pmatrix}
3\\t_{m-2}\end{pmatrix}-\begin{pmatrix}
3\\t_{m-1}\end{pmatrix}-\begin{pmatrix}
3\\t_m\end{pmatrix}+\begin{pmatrix}
3\\t_{m-2}\end{pmatrix}\begin{pmatrix}
3\\t_{m-1}\end{pmatrix}\right.+
$$
$$
\left.+\begin{pmatrix} 3\\t_{m-2}\end{pmatrix}\begin{pmatrix}
3\\t_m\end{pmatrix}+\begin{pmatrix}
3\\t_{m-1}\end{pmatrix}\begin{pmatrix}
3\\t_m\end{pmatrix}-\begin{pmatrix}
3\\t_{m-2}\end{pmatrix}\begin{pmatrix}
3\\t_{m-1}\end{pmatrix}\begin{pmatrix} 3\\t_m\end{pmatrix}\right)=
$$
$$
=\begin{cases} (-1)^m,\;\; if \;\; t_m\geq 4\\
0, \;\; if \;\; t_m=3\\
0, \;\; if \;\; t_m=2,\;\;t_{m-1}=3\\
2(-1)^{m-1}, \;\; if \;\; t_m=2,\;\;t_{m-1}> 3\\
0, \;\; if \;\; t_m=1,\;\;t_{m-1}=2,\;\;t_{m-2}=3\\
4(-1)^m, \;\; if \;\; t_m=1,\;\;t_{m-1}=2,\;\ t_{m-2}> 3\\
0, \;\; if \;\; t_m=1,\;\;t_{m-1}=3,\;\;t_{m-2}>3\\
2(-1)^{m-1}, \;\; if \;\; t_m=1,\;\;t_{m-1}> 3\end{cases} ;
$$
etc.\newline
\indent We see that, for every fixed $n,$ we have a finite number of distinct values of the sequence $\{\left\{\begin{matrix} n\\
k\end{matrix}\right\}\}_{k\geq0}.$ Therefore it is bounded by a constant $C(n).\blacksquare$
Denote, for any $n\in\mathbb{N},$

\begin{equation}\label{56}
F(n,x)=\sum^\infty_{k=0}\left \{\begin{matrix} n\\k
\end{matrix}\right \}x^k,\;\; |x|<1.
\end{equation}

Put
\begin{equation}\label{57}
\tau(x)=\sum^\infty_{k=0}\tau_k x^k,\;\; |x|<1,
\end{equation}

where ${\tau}$ is the Thue-Morse sequence (\ref{16}).

\begin{theorem}\label{t18}    For every $n\in\mathbb{N}$ the
quotient $\frac{F(n,x)}{\tau(x)}$ is a rational function.
\end{theorem}

\bfseries Proof. \mdseries  It follows from Theorems \ref{t9} and
\ref{t14} that

\begin{equation}\label{58}
\left \{\begin{matrix} n\\k
\end{matrix}\right \}=\tau_k+\sum_{1\leq i\leq\log_2(2k)}\begin{pmatrix} n\\ i \end{pmatrix}
\left \{\begin{matrix} i\\k-2^{i-1}
\end{matrix}\right \},\;\; k\geq 1.
\end{equation}

Note that, for $k=0,$ we have

$$
\left \{\begin{matrix} n\\0
\end{matrix}\right \}=1=\tau_0.
$$

Therefore, by (\ref{56})-(\ref{58}), we find\newpage

$$
F(n,x)=\sum^\infty_{k=0}\left \{\begin{matrix} n\\k
\end{matrix}\right \}x^k=\tau(x)+\sum_{k\geq 1}\sum_{1\leq i \leq\log_2(2k)}
\begin{pmatrix} n\\ i \end{pmatrix}\left \{\begin{matrix} i\\k-2^{i-1}
\end{matrix}\right \}x^k=
$$
$$
=\tau(x)+\sum^n_{i=1}\begin{pmatrix} n\\ i
\end{pmatrix}\sum^\infty_{k=2^{i-1}}x^k\left \{\begin{matrix}
i\\k-2^{i-1}
\end{matrix}\right \}=
$$
$$
=\tau(x)+\sum^n_{i=1}\begin{pmatrix} n\\ i
\end{pmatrix}\sum^\infty_{r=0}\left \{\begin{matrix}
i\\r
\end{matrix}\right \}x^{r+2^{i-1}}=
$$
\begin{equation}\label{59}
=\tau(x)+\sum^n_{i=1}\begin{pmatrix} n\\i
\end{pmatrix}x^{2^{i-1}}F(i,x),\;\;\;|x|< 1.
\end{equation}

Thus we obtain us a recursion formula for $F(n,x):$

\begin{equation}\label{60}
\left(1-x^{2^{n-1}}\right)
F(n,x)=\tau(x)+\sum^{n-1}_{i=1}\begin{pmatrix} n\\i
\end{pmatrix}x^{2^{i-1}}F(i,x),\;\; |x|< 1.
\end{equation}

Put

\begin{equation}\label{61}
F(n,x)=\tau(x)\frac{P_n(x)}{(1-x)(1-x^2)\hdots
(1-x^{2^{n-2}})}\;\;(|x|< 1).
\end{equation}

Then we obtain a recursion formula for $P_n(x):$

$$
P_n(x)=\frac{1}{1-x^{2^{n-1}}}\left((1-x)(1-x^2)\hdots
(1-x^{2^{n-1}})+\right.
$$
\begin{equation}\label{62}
\left.+\sum^{n-1}_{i=1}\begin{pmatrix} n\\ i
\end{pmatrix}(1-x^{2^i})(1-x^{2^{i+1}})\hdots
(1-x^{2^{n-1}})\;x^{2^{i-1}}P_i(x)\right).
\end{equation}

Here we do not cancel $1-x^{2^{n-1}}$ in order to avoid
additional conventions. In particular, from (62) we find

$$
P_1(x)=\frac{1}{1-x}(1-x)=1,
$$

$$
P_2(x)=\frac{1}{1-x^2}((1-x)(1-x^2)+2(1-x^2)x)=1+x,
$$

$$
P_3(x)=1+2x+2x^2+x^3,
$$

$$
P_4(x)=1+3x+5x^2+3x^3+3x^4+5x^5+3x^6+x^7,
$$
\newpage

$$
P_5(x)=1+4x+9x^2+6x^3+9x^4+16x^5+11x^6+4x^7+
$$

$$
4x^8+11x^9+16x^{10}+9x^{11}+6x^{12}+9x^{13}+4x^{14}+x^{15},
$$

$$
P_6(x)=1+5x+14x^2+10x^3+19x^4+35x^5+26x^6+10x^7+14x^8+40x^9+
$$

$$
+61x^{10}+35x^{11}+26x^{12}+40x^{13}+19x^{14}+5x^{15}+5x^{16}+19x^{17}+40x^{18}+
$$

$$
+26x^{19}+35x^{20}+61x^{21}+40x^{22}+14x^{23}+10x^{24}+26x^{25}+35x^{26}+19x^{27}+
$$

$$
+10x^{28}+14x^{29}+5x^{30}+x^{31},
$$
etc.

By simple induction we see that $P_n(x)$ is a polynomial in $x$ of
degree $2^{n-1}-1$. Thus, the theorem follows from
(\ref{61}).$\blacksquare$

     But (\ref{61}) gives us more in view of identity:

$$
(1-x)(1-x^2)\hdots (1-x^{2^{n-2}})=\sum^{2^{n-1}-1}_{k=0}\tau_k x^k=
\tau(x) + o(x^{2^{n-1}-1}).
$$

Therefore, from (\ref{61}) it follows that

$$
P_n(x)=F(n,x)(1+o(x^{2^{n-1}-1})).
$$

Thus polynomial $P_n(x)$ of degree $2^{n-1}-1$ is a partial sum of (\ref{56}). Therefore, we obtain the following result.

\begin{theorem}\label{t19}
Polynomial $P_n(x)$ which is defined recursively by $(\ref{62})$ is
equal to

\begin{equation}\label{63}
P_n(x)=\sum^{2^{n-1}-1}_{k=0}\left\{\begin{matrix} n\\k
\end{matrix}\right\} x^k.
\end{equation}
\end{theorem}
\indent Moreover, from (\ref{61}),for every $n\in \mathbb{N},$ we have an identity
\newpage

\begin{equation}\label{64}
\frac{\sum\limits^\infty\limits_{k=0}\left\{\begin{matrix} n\\k
\end{matrix}\right\}x^k}{\sum\limits^{2^{n-1}-1}\limits_{k=0}\left\{\begin{matrix} n\\k
\end{matrix}\right\}x^k}=\prod\limits^\infty\limits_{i=n-1}(1-x^{2^i})\;\;(|x|< 1).
\end{equation}

\section{Another type of recursion for basis polynomials}

For any $k\in\mathbb{N},$ let us consider the set $A_k$ of those
positive integers $i\leq \log_2(2k)$ for which
$\left\lfloor\frac{k}{2^i}-\frac 1 2\right\rfloor=\left\lceil
\frac{k+1}{2^i}-1\right\rceil.$
The common values of these expressions denote by
$\lambda(k;i):$
\begin{equation}\label{65}
\lambda(k;i)=\left\lfloor\frac{k}{2^i}-\frac 1 2\right\rfloor=\left\lceil
\frac{k+1}{2^i}-1\right\rceil.
\end{equation}

\begin{theorem}\label{t20}
\begin{equation}\label{66}
\left\{\begin{matrix} n\\k
\end{matrix}\right\}=\tau_k+\sum\limits_{i\in A_k}\begin{pmatrix}n\\ i \end{pmatrix}
\left\{\begin{matrix} i\\k-2^{i-1}-\lambda(k;i)2^i
\end{matrix}\right\}\tau_{\lambda(k;,i)}.
\end{equation}
\end{theorem}

\bfseries Proof. \mdseries   Taking into account (\ref{63}) and
comparing $Coef_{x^k},\enskip k\leq 2^{n-1}-1,$ in both sides of (\ref{62})
we find

$$
\left\{\begin{matrix} n\\k
\end{matrix}\right\}=\tau_k+\sum\limits^{n-1}\limits_{i=1}\begin{pmatrix}n\\ i \end{pmatrix}
\sum_l \tau_l \left\{\begin{matrix} i\\k-2^{i-1}-2^i l
\end{matrix}\right\},
$$

where the summing is over those values of $l\geq 0$ for which

$$
2^{i-1}+2^i l\leq k,\;\;k-2^{i-1}-2^i l \leq deg P_i(x)=2^{i-1}-1.
$$

Consequently,

$$
l\in\left[\frac{k+1}{2^i}-1,\;\;\frac{k}{2^i}-\frac 1 2 \right],
$$

Nevertheless, the length of this segment equals to $\frac 1 2
-\frac{1}{2^i}$. This means that
$$l=\lambda(k;i)= \left\lfloor\frac{k}{2^i}-\frac 1 2
\right\rfloor=\left\lceil\frac{k+1}{2^i}-1\right\rceil\geq 0.\blacksquare$$

\begin{example}

Let $k=2^m$. Then $i\leq \log_2(2k)=m+1$. If $i\leq m$ then
$\left\lfloor\frac{k}{2^i}-\frac 1 2\right\rfloor=2^{m-i}-1$ while
$\left\lceil\frac{k+1}{2^i}-1\right\rceil=2^{m-i}$. It is left to
consider the case $i=m+1$ for which

$$
\left\lfloor\frac{k}{2^i}-\frac 1 2
\right\rfloor=\left\lceil\frac{k+1}{2^i}-1\right\rceil=0.
$$

Thus, by $(\ref{66}),$ we have\newpage

$$
\left\{\begin{matrix} n\\ 2^m\end{matrix}\right\}=-1+\begin{pmatrix}
n\\ m+1 \end{pmatrix}\left\{\begin{matrix} m+1\\
0\end{matrix}\right\}=\begin{pmatrix} n\\ m+1 \end{pmatrix}-1.
$$
\end{example}

\section{Conversion of Theorem 12}

We have seen (Theorem 12) that

 \begin{equation}\label{67}
 If\enskip b\equiv 2^{a-1}+c\pmod{2^a},\;\;c\in[0,2^{a-1}),\enskip
then\enskip \left\{\begin{matrix} a\\ b\end{matrix}\right\}=0.
\end{equation}

The conversion of Theorem 12 is based on the following result.

\begin{theorem}\label{t21}
Let $a\in\mathbb{N},\;\; r\geq a,\;\;c\in[0,2^{a-1})$. Then for
$b=2^rl+c,$ where $l$ is odd, we have

\begin{equation}\label{68}
\left\{\begin{matrix} a\\
b\end{matrix}\right\}=\tau_l\left\{\begin{matrix} a\\
c\end{matrix}\right\}.
\end{equation}
\end{theorem}

\bfseries Proof. \mdseries Let
$k=2^{t_1-1}+2^{t_2-1}+\hdots+2^{t_m-1},\;\;t_1>t_2>\hdots >t_m\geq
1$. Comparing Theorems 9 and 20, we conclude that

\begin{equation}\label{69}
 A_k=\{t_1,t_2,\hdots,t_m\}
\end{equation}

and, for $i\in A_k,$

\begin{equation}\label{70}
\left\{\begin{matrix} i\\
k-2^{i-1}-\lambda(k;i)2^i\end{matrix}\right\}\tau_{\lambda(k;i)}=
\left\{\begin{matrix} i\\
k-2^{i-1}\end{matrix}\right\}.
\end{equation}

In particular, for $i=t_j$ we have

$$
\lambda(k;t_j)=\left\lfloor\frac{k}{2^{t_j}}-\frac 1 2
\right\rfloor=\left\lceil\frac{k+1}{2^{t_j}}-1\right\rceil.
$$

Thus,

$$
\lambda(k;t_j)=\left\lfloor\frac{2^{t_1-1}+2^{t_2-1}+\hdots+2^{t_j-1}+\hdots+
2^{t_m-1}}{2^{t_j}}-\frac 1 2\right\rfloor=
$$
\begin{equation}\label{71}
=2^{t_1-t_j-1}+2^{t_2-t_j-1}+\hdots+2^{t_{j-1}-t_j-1}
\end{equation}

and, consequently,

\begin{equation}\label{72}
\tau_{\lambda(k;\enskip t_j)}=(-1)^{j-1}.
\end{equation}

Therefore, by (\ref{70}) (for $i=t_j$) and (\ref{71}), we find

\begin{equation}\label{73}
(-1)^{j-1}\left\{\begin{matrix} t_j\\
k-2^{t_1-1}-2^{t_2-1}-\hdots-2^{t_{j-1}-1}-2^{t_l-1}\end{matrix}\right\}=
\left\{\begin{matrix} t_j\\
k-2^{t_j-1}\end{matrix}\right\}
\end{equation}\newpage

or, taking into account that $k=2^{t_1-1}+\hdots+2^{t_m-1},$

$$(-1)^{j-1}\left\{\begin{matrix} t_j\\
2^{t_{j+1}-1}+\hdots+2^{t_m-1}\end{matrix}\right\}=$$
\begin{equation}\label{74}
\left\{\begin{matrix} t_j\\
2^{t_1-1}+\hdots+2^{t_{j-1}-1}+2^{t_{j+1}-1}+\hdots+2^{t_m-1}\end{matrix}\right\}.
\end{equation}

Put here

$$
t_j=a,\;\;2^{t_{j+1}-1}+\hdots+2^{t_m-1}=c\in[0, 2^{a-1}),
$$

$$
t_{j-1}-1=r\geq
a,\;\;2^{t_1-1}+\hdots+2^{t_{j-1}-1}=2^rl,
$$
where $$l=2^{t_1-1-r}+2^{t_2-1-r}+...+2^{t_{j-1}-1-r}=2^{t_1-1-r}+2^{t_2-1-r}+...+1.$$
In view of $\tau_l=(-1)^{j-1},$
and
$$
b=2^rl+c,
$$
we write (\ref{74}) in the form of (\ref{68}) $\blacksquare$

Since in Theorem 21 $b-c=2^rl$ then $\tau_l=\tau_{b-c}$. Therefore,
Theorem 21, for $b=k,$ one can write in the following form.

\begin{corollary}\label{22}
$$
\left\{\begin{matrix} a\\
k\end{matrix}\right\}=\tau_{k-i}\left\{\begin{matrix} a\\
i\end{matrix}\right\},\;\; k\equiv
i\pmod{2^a},\;\;i=0,1,2,\hdots,2^{a-1}-1.
$$
\end{corollary}
It is worth to add that by (\ref{15})

\begin{equation}\label{75}
\left\{\begin{matrix} 0\\
k\end{matrix}\right\}=\tau_k.
\end{equation}

In particular, taking into account (\ref{67}) and (\ref{75}), we
obtain the following sequences:

$$
\left\{\begin{matrix} 0\\
k\end{matrix}\right\}: 1,\;-1,\;-1,\;1,\;-1,\;1,\;1,\;-1,\;\hdots
$$
$$
\left\{\begin{matrix} 1\\
k\end{matrix}\right\}:\mathbf
1,0,-1,0,-1,0,1,0,-1,0,1,0,1,0,-1,\hdots
$$
$$
\left\{\begin{matrix} 2\\
k\end{matrix}\right\}:\mathbf
1,\mathbf1,0,0,-1,-1,0,0,-1,-1,0,0,1,1,\hdots
$$
$$
\left\{\begin{matrix} 3\\
k\end{matrix}\right\}:\mathbf 1,\mathbf 2,\mathbf 2,\mathbf
1,0,0,0,0,-1,-2,-2,-1,0,0,0,0,\hdots
$$\newpage
$$
\left\{\begin{matrix} 4\\
k\end{matrix}\right\}:\mathbf {1, 3,5, 3, 3, 5, 3,
1},\underbrace{0,\hdots,0}_8,-1,-3,-5,-3,-3,\hdots
$$
$$
\left\{\begin{matrix} 5\\
k\end{matrix}\right\}\mathbf
{1,4,9,6,9,16,11,4,4,11,16,9,6,9,4,1,}\underbrace{0,\hdots,0}_{16},-1,\hdots
$$
$$
\left\{\begin{matrix} 6\\
k\end{matrix}\right\}\mathbf{1,5,14,10,19,35,26,10,14,40,61,35,26,40,19,5,5,19,40,26,35,61,}
$$
$$
\;\;\;\;\;\mathbf{40,14,10,26,35,19,10,14,5,1,}\underbrace{0,\hdots,0}_{32},-1,-5,-14,-10,-19,-35,\hdots
$$
etc.\newline
\indent We see that, if in sequence $(\left\{\begin{matrix} a\\
k\end{matrix}\right\}) $ to ignore the signs, then it becomes to periodic sequence with period $2^a$ such that the second part of period containing terms of the form $2^{a-1}+c,\enskip c=1,...,2^{a-1},$ by (67), consists of 0's. It is left to show that all terms $\left\{\begin{matrix} a\\
k\end{matrix}\right\},\enskip k=1,...,2^{a-1},$ of the first part of period are positive. Note that, in order to prove this, it is sufficient to find at least one permutation of elements $1,...,a$ with a given index $k\in[1,2^{a-1}].$ We can consider only positions of 0's in $a-1$-digit binary expansion of $k.$ Start with initial permutation $\pi=(1,2,...,a).$
Suppose that the first (from the left) run of binary 0's of $k$ has positions $i,...,i+t-1.$ Let us write in $\pi$ the elements $i,...,i+t,i+t$ in the reverse order. Then we obtain permutation $\pi_1=(1,...,i-1,i+t,i+t-1,...,i+1,i,i+t+1,...,n),$ having the first (from the right) run of descent points exactly in the required positions: $\pi_1(i)=i+t>\pi_1(i+1)=i+t-1>...>\pi_1(i+t-1)=i+1>\pi_(i+t)=i.$ Further, by the same algorithm, we obtain permutation with the first two runs of descent points corresponding to the first two runs of 0's in the binary expansion of $k,$ etc. At the end of this process we obtain a permutation with index $k.$ So, we proved that $\left\{\begin{matrix} a\\k\end{matrix}\right\}\geq1,\enskip k=1,...,2^{a-1}.$ This gives conversion of Theorem 13 and, consequently, of Theorem 12.
Thus we have the following statement.

\begin{theorem}\label{t22}
If $k=2^{t_1-1}+2^{t_2-1}+\hdots+2^{t_m-1},\;t_1>t_2>\hdots>t_m\geq
1$, then $t_i,\; i=1,2,\hdots, m$, are only positive integer roots
of the polynomials $\left\{\begin{matrix} n\\
k\end{matrix}\right\}$.
\end{theorem}

\begin{corollary}\label{t2} Let $k=2^{t_1-1}+...+2^{t_m-1},\enskip t_1>...>t_m.$ If $n\neq t_i,\enskip i=1,...,m,$ then\newpage
$$ rank{ \begin{pmatrix} 1& \begin{pmatrix}n\\t_1\end{pmatrix}&
\begin{pmatrix}n\\t_2\end{pmatrix}& \begin{pmatrix}n\\t_3\end{pmatrix}& \hdots &
\begin{pmatrix}n\\t_m\end{pmatrix}\\
1& 1 &\begin{pmatrix}t_1\\t_2\end{pmatrix}&
\begin{pmatrix}t_1\\t_3\end{pmatrix}&
\hdots &\begin{pmatrix}t_1\\t_m\end{pmatrix}\\
1& 0& 1 &\begin{pmatrix}t_2\\t_3\end{pmatrix}&\hdots
&\begin{pmatrix}t_2\\t_m\end{pmatrix}\\1 & 0 & 0 & 1 & \hdots &
\begin{pmatrix}t_3\\t_m\end{pmatrix}\\\hdots &\hdots &\hdots &\hdots &\hdots &\hdots
\\1 & 0 & 0 & 0 & \hdots & 1
\end{pmatrix} }=m+1.$$
\end{corollary}
\bfseries Proof\mdseries \enskip follows from Theorem \ref{t22} and determinant representation (\ref{40}) of basis polynomials. $\blacksquare$ \newline
\indent In conclusion, note that the bold-faced numbers are coefficients of the
corresponding polynomials $P_i(x),\;\;i=1,2,\hdots $ which are
defined by recursion (\ref{62}).

If consider the concatenation sequence of the all bold-faced numbers
(cf.\cite{15},A060351):

$$
\left\{\begin{matrix} 1\\
0\end{matrix}\right\};\left\{\begin{matrix} 2\\
0\end{matrix}\right\},\left\{\begin{matrix} 2\\
1\end{matrix}\right\};\left\{\begin{matrix} 3\\
0\end{matrix}\right\},\left\{\begin{matrix} 3\\
1\end{matrix}\right\},\left\{\begin{matrix} 3\\
2\end{matrix}\right\},\left\{\begin{matrix} 3\\
3\end{matrix}\right\};\left\{\begin{matrix} 4\\
0\end{matrix}\right\},\left\{\begin{matrix} 4\\
1\end{matrix}\right\},\left\{\begin{matrix} 4\\
2\end{matrix}\right\},
$$
\begin{equation}\label{76}
\left\{\begin{matrix} 4\\
3\end{matrix}\right\},\left\{\begin{matrix} 4\\
4\end{matrix}\right\},\left\{\begin{matrix} 4\\
5\end{matrix}\right\},\left\{\begin{matrix} 4\\
6\end{matrix}\right\},\left\{\begin{matrix} 4\\
7\end{matrix}\right\};\hdots
\end{equation}
then it is easy to see that this sequence one can write in an explicit form:

\begin{equation}\label{77}
\{ \left\{\begin{matrix} \lfloor\log_2{k}\rfloor+1\\
k-2^{\lfloor\log_2{k}\rfloor}\end{matrix}\right\} \}^\infty_{k=1}.
\end{equation}

This sequence is closely connected with asymptotics of $\left\{\begin{matrix} n\\
k\end{matrix}\right\}$.

\begin{theorem}\label{t23}
$$
\left\{\begin{matrix} n\\
k\end{matrix}\right\}\sim \left\{\begin{matrix} \lfloor\log_2{k}\rfloor+1\\
k-2^{\lfloor\log_2{k}\rfloor}\end{matrix}\right\}\begin{pmatrix}
n\\\lfloor\log_2{k}\rfloor+1 \end{pmatrix}\sim
$$

\begin{equation}\label{78}
\sim\frac{1}{\lfloor\log_2{k}+1\rfloor!}\{\left\{\begin{matrix} \lfloor\log_2{k}\rfloor+1\\
k-2^{\lfloor\log_2{k}\rfloor}\end{matrix}\right\}n^{\lfloor\log_2{k}\rfloor+1}
\;\;(n\rightarrow\infty)
\end{equation}
\end{theorem}

\bfseries Proof.\mdseries The theorem follows directly from
(\ref{58})$\blacksquare$

Thus, according to (\ref{78}), the first coefficients of $\left\{\begin{matrix} n\\
k\end{matrix}\right\},\;k\geq 1$, as linear combinations of binomial
coefficients form sequence (\ref{76}) ((cf. Appendix) ).

\newpage
\begin{remark}\label{4}  Some observations on the basis polynomials have been
done by the author as early as 1993 \cite{13} with the calculation
of some first polynomials. But only in the current paper we give a
foundation of a more perfect theory of these polynomials.
\end{remark}

\section{Some open problems}

1. We conjecture that all real roots of the basis polynomials are
rational.

2. We conjecture that a polynomial $\left\{\begin{matrix} n\\
k\end{matrix}\right\},\; k\geq 1$, has only real roots if and only
if the number of 0's in the binary expansion of $k$ less that $2$.
In view of Theorem 12 this condition is sufficient (since there is
no a place for two conjugate complex roots). Therefore, it is left
to prove its necessity. We verified this conjecture up to
$k=32$. We have polynomials $\left\{\begin{matrix} n\\
k\end{matrix}\right\}$ with only real roots for

$$
k=1,2,3,5,6,7,11,13,14,15,23,27,29,30,31,\hdots
$$
(cf. sequence A 089633 \cite{15}).

3. It is interesting to investigate the sequence $\{k_m\}$ for which
the polynomials $\left\{\begin{matrix} n\\
k_m\end{matrix}\right\}$ have a root $n=-1$. The first values of
$k_m$ are: $2,5,8,11,23,\hdots$

4. Let $D_n^{(a)}$ be the number of alternating permutations without
fixed points (i.e.$\pi(i)\neq i,\;\;i=1,2,\hdots,n $). We conjecture
that
$$
\lim_{n\rightarrow \infty}\frac{D_n^{(a)}}{a_n}=e^{-1},
$$
where ${a_n}$ is the sequence (\ref{21}).

5. Let $S^{(a)}(n,l)$ be the number of alternating permutations
having $l$ cycles (the absolute value of the "alternating" Stirling
numbers of the first kind). We conjecture that for a fixed $l$

$$
\lim_{n\rightarrow\infty}\frac{nS^{(a)}(n,l)}{a_n(\ln{n})^{l-1}}=\frac{1}{(l-1)!},
$$

where $a_n$  is the sequence (\ref{21}). The latter means that for
each $l$ the events "a permutation is alternative" and "a
permutation has $l$ cycles" are asymptotically independent.

\; \;

\bfseries Appendix. \mdseries
List of the first 32 basis polynomials
$$
\left\{\begin{array}{c} n \\ 0 \end{array}\right\}=1
$$
$$
\left\{\begin{array}{c} n \\ 1
\end{array}\right\}=\left(\begin{array}{c} n \\ 1
\end{array}\right)-1
$$
$$
\left\{\begin{array}{c} n \\ 2
\end{array}\right\}=\left(\begin{array}{c} n \\ 2
\end{array}\right)-1
$$
$$
\left\{\begin{array}{c} n \\ 3
\end{array}\right\}=\left(\begin{array}{c} n \\ 2
\end{array}\right)-\left(\begin{array}{c} n \\ 1
\end{array}\right)+1
$$
$$
\left\{\begin{array}{c} n \\ 4
\end{array}\right\}=\left(\begin{array}{c} n \\ 3
\end{array}\right)-1
$$
$$
\left\{\begin{array}{c} n \\ 5
\end{array}\right\}=2\left(\begin{array}{c} n \\ 3
\end{array}\right)-\left(\begin{array}{c} n \\ 1
\end{array}\right)+1
$$\newpage
$$
\left\{\begin{array}{c} n \\ 6
\end{array}\right\}=2\left(\begin{array}{c} n \\ 3
\end{array}\right)-\left(\begin{array}{c} n \\ 2
\end{array}\right)+1
$$
$$
\left\{\begin{array}{c} n \\ 7
\end{array}\right\}=\left(\begin{array}{c} n \\ 3
\end{array}\right)-\left(\begin{array}{c} n \\ 2
\end{array}\right)+\left(\begin{array}{c} n \\ 1
\end{array}\right)-1
$$
$$
\left\{\begin{array}{c} n \\ 8
\end{array}\right\}=\left(\begin{array}{c} n \\ 4
\end{array}\right)-1
$$
$$
\left\{\begin{array}{c} n \\ 9
\end{array}\right\}=3\left(\begin{array}{c} n \\ 4
\end{array}\right)-\left(\begin{array}{c} n \\ 1
\end{array}\right)+1
$$
$$
\left\{\begin{array}{c} n \\ 10
\end{array}\right\}=5\left(\begin{array}{c} n \\ 4
\end{array}\right)-\left(\begin{array}{c} n \\ 2
\end{array}\right)+1
$$
$$
\left\{\begin{array}{c} n \\ 11
\end{array}\right\}=3\left(\begin{array}{c} n \\ 4
\end{array}\right)-\left(\begin{array}{c} n \\ 2
\end{array}\right)+\left(\begin{array}{c} n \\ 1
\end{array}\right)-1
$$
$$
\left\{\begin{array}{c} n \\ 12
\end{array}\right\}=3\left(\begin{array}{c} n \\ 4
\end{array}\right)-\left(\begin{array}{c} n \\ 3
\end{array}\right)+1
$$
$$
\left\{\begin{array}{c} n \\ 13
\end{array}\right\}=5\left(\begin{array}{c} n \\ 4
\end{array}\right)-2\left(\begin{array}{c} n \\ 3
\end{array}\right)+\left(\begin{array}{c} n \\ 1
\end{array}\right)-1
$$
$$
\left\{\begin{array}{c} n \\ 14
\end{array}\right\}=3\left(\begin{array}{c} n \\ 4
\end{array}\right)-2\left(\begin{array}{c} n \\ 3
\end{array}\right)+\left(\begin{array}{c} n \\ 2
\end{array}\right)-1
$$
$$
\left\{\begin{array}{c} n \\ 15
\end{array}\right\}=\left(\begin{array}{c} n \\ 4
\end{array}\right)-\left(\begin{array}{c} n \\ 3
\end{array}\right)+\left(\begin{array}{c} n \\ 2
\end{array}\right)-\left(\begin{array}{c} n \\ 1
\end{array}\right)+1
$$
$$
\left\{\begin{array}{c} n \\ 16
\end{array}\right\}=\left(\begin{array}{c} n \\ 5
\end{array}\right)-1
$$
$$
\left\{\begin{array}{c} n \\ 17
\end{array}\right\}=4\left(\begin{array}{c} n \\ 5
\end{array}\right)-\left(\begin{array}{c} n \\ 1
\end{array}\right)+1
$$
$$
\left\{\begin{array}{c} n \\ 18
\end{array}\right\}=9\left(\begin{array}{c} n \\ 5
\end{array}\right)-\left(\begin{array}{c} n \\ 2
\end{array}\right)+1
$$
$$
\left\{\begin{array}{c} n \\ 19
\end{array}\right\}=6\left(\begin{array}{c} n \\ 5
\end{array}\right)-\left(\begin{array}{c} n \\ 2
\end{array}\right)+\left(\begin{array}{c} n \\ 1
\end{array}\right)-1
$$
$$
\left\{\begin{array}{c} n \\ 20
\end{array}\right\}=9\left(\begin{array}{c} n \\ 5
\end{array}\right)-\left(\begin{array}{c} n \\ 3
\end{array}\right)+1
$$
$$
\left\{\begin{array}{c} n \\ 21
\end{array}\right\}=16\left(\begin{array}{c} n \\ 5
\end{array}\right)-2\left(\begin{array}{c} n \\ 3
\end{array}\right)+\left(\begin{array}{c} n \\ 1
\end{array}\right)-1
$$
$$
\left\{\begin{array}{c} n \\ 22
\end{array}\right\}=11\left(\begin{array}{c} n \\ 5
\end{array}\right)-2\left(\begin{array}{c} n \\ 3
\end{array}\right)+\left(\begin{array}{c} n \\ 2
\end{array}\right)-1
$$
$$
\left\{\begin{array}{c} n \\ 23
\end{array}\right\}=4\left(\begin{array}{c} n \\ 5
\end{array}\right)-\left(\begin{array}{c} n \\ 3
\end{array}\right)+\left(\begin{array}{c} n \\ 2
\end{array}\right)-\left(\begin{array}{c} n \\ 1
\end{array}\right)+1
$$
$$
\left\{\begin{array}{c} n \\ 24
\end{array}\right\}=4\left(\begin{array}{c} n \\ 5
\end{array}\right)-\left(\begin{array}{c} n \\ 4\end{array}\right)
+1
$$
$$
\left\{\begin{array}{c} n \\ 25
\end{array}\right\}=11\left(\begin{array}{c} n \\ 5
\end{array}\right)-3\left(\begin{array}{c} n \\ 4
\end{array}\right)+\left(\begin{array}{c} n \\ 1
\end{array}\right)-1
$$\newpage
$$
\left\{\begin{array}{c} n \\ 26
\end{array}\right\}=16\left(\begin{array}{c} n \\ 5
\end{array}\right)-5\left(\begin{array}{c} n \\ 4
\end{array}\right)+\left(\begin{array}{c} n \\ 2
\end{array}\right)-1
$$
$$
\left\{\begin{array}{c} n \\ 27
\end{array}\right\}=9\left(\begin{array}{c} n \\ 5
\end{array}\right)-3\left(\begin{array}{c} n \\ 4
\end{array}\right)+\left(\begin{array}{c} n \\ 2
\end{array}\right)-\left(\begin{array}{c} n \\ 1
\end{array}\right)+1
$$
$$
\left\{\begin{array}{c} n \\ 28
\end{array}\right\}=6\left(\begin{array}{c} n \\ 5
\end{array}\right)-3\left(\begin{array}{c} n \\ 4
\end{array}\right)+\left(\begin{array}{c} n \\ 3
\end{array}\right)-1
$$
$$
\left\{\begin{array}{c} n \\ 29
\end{array}\right\}=9\left(\begin{array}{c} n \\ 5
\end{array}\right)-5\left(\begin{array}{c} n \\ 4
\end{array}\right)+2\left(\begin{array}{c} n \\ 3
\end{array}\right)-\left(\begin{array}{c} n \\ 1
\end{array}\right)+1
$$
$$
\left\{\begin{array}{c} n \\ 30
\end{array}\right\}=4\left(\begin{array}{c} n \\ 5
\end{array}\right)-3\left(\begin{array}{c} n \\ 4
\end{array}\right)+2\left(\begin{array}{c} n \\ 3
\end{array}\right)-\left(\begin{array}{c} n \\ 2
\end{array}\right)+1
$$
$$
\left\{\begin{array}{c} n \\ 31
\end{array}\right\}=\left(\begin{array}{c} n \\ 5
\end{array}\right)-\left(\begin{array}{c} n \\ 4
\end{array}\right)+\left(\begin{array}{c} n \\ 3
\end{array}\right)-\left(\begin{array}{c} n \\ 2
\end{array}\right)+\left(\begin{array}{c} n \\ 1
\end{array}\right)-1
$$

\begin{thebibliography}{18}
\bibitem 1. M.Abramowitz and I.A.Stegun (Eds.), Bernoulli and Euler Polynomials and the
Euler-Maklaurin Formula in \slshape Handbook of Manhematical
Functions with Formulas, Graphs, and Manhematical Tables, 9th
printing.\upshape ,New York: Dover,pp. 804-806, 1972.\newpage
\bibitem 2. D.Andre, Sur les Permutation Alternees,  \slshape J.Math.
Pures Appl., \upshape  \bfseries 7 \mdseries (1881),167-184.
\bibitem 3.V.Arnold, Bernoulli-Euler updown numbers associated with
function singularities, their combinatorics and arithmetic,\slshape
Duke Math.J.   \upshape \bfseries 63(2) \mdseries (1990), 537-555.
\bibitem 4. F.C.S.Brown, T.M.A. Fink, K.Willbrand, On arithmetic and
asymptotic properties of up-down numbers, \slshape Discrete
Math.\upshape \bfseries 307 \mdseries (2007),1722-1736.
\bibitem 5. N.G.de Bruijn, Permutations with given ups and downs,
\slshape Nieuw Arch.\upshape \bfseries  3 \mdseries(1970), 61-65.
\bibitem 6. L.Carlitz, Permutations with prescribed pattern, II,
\slshape Matem. Nachr.\upshape \bfseries 83 \mdseries (1978),
101-126.
\bibitem 7. H.O.Foulkes, Enumeration of permutations with prescribed up-down and
 inversion sequences,\slshape Discrete Math.\upshape \bfseries 15 \mdseries (1976),
235-252.
\bibitem 8. S.Goldstein, K.A.Kelly and E.R.Speer,The fractal structure of rarefied
sums of the Thue-Morse sequence, \slshape J.Number Th.\upshape
\bfseries 42 \mdseries (1992), 1-19.
\bibitem 9. C.L.Mallows, L.A.Shepp, Enumerating pairs of permutations with the same
up-down form, \slshape Descrete Math. \upshape \bfseries
54\mdseries(1985), 301-311.
\bibitem {10}. M.Morse, Reccurent geodesics on a surface of negative
curvature,\slshape  Trans. Amer.Math.Soc.\upshape \bfseries 22
\mdseries(1921), 84-100.
\bibitem {11}. I.Niven, A combinatorial problem of finite sequences,
\slshape Nieuw Arch.Wisk \upshape \bfseries 3 \mdseries (1968),
116-123.
\bibitem {12}. B.Shapiro, M.Shapiro, A.Vainshtein, Periodic de Bruijn
triangles: exact and asymptotic results, \slshape Discrete Math.
\upshape \bfseries 298(1-3) \mdseries (2005), 321-333.
\bibitem {13}. V.S.Shevelev, A classification of permutations by its geometric structure,
\slshape Deposed in VINITI, no.1457-B93 (1993), 1-20 (in Russian).
\bibitem {14}. V.S.Shevelev, On an arithmetic property of permutation numbers with a given
signature associated with the Morse sequence,\slshape
Izv.Vyssh.Uchebn.Zaved. Sev.-Kavk.Reg.Estestv.Nauki \upshape
\bfseries 2 \mdseries (1996), 20-24 (in Russian; MR99e: 11023).


\bibitem {15}.N.J.A.Sloane, The On-Line Encyclopedia of Integer
Sequences (http://www.research.att.com:/$\sim$njas /sequenses/).
\bibitem {16}. R. Stanley, Alternating permutations and symmetric functions, J. Combin. Theory Series A \bfseries 114 \mdseries (2007), 436-460.
\bibitem {17}. G.Szpiro, The number of permutations with a given
signature and the expectations of their elements,\slshape Discrete
Math.\upshape \bfseries 226 \mdseries (2001), 423-430.
\bibitem {18}. G.Viennot. Permutations ayant une forme donnee.\slshape Discrete
Math.\upshape \bfseries 26 \mdseries (1979), 279-284.
\end{thebibliography}
\end {document}